\documentclass[11pt]{amsart}

\usepackage{amsmath, amssymb, amsfonts, amsthm}
\usepackage{graphicx}
\usepackage{color}
\usepackage{bm}
\usepackage{subfigure}
\usepackage{graphicx}
\usepackage{mathabx}
\usepackage{multirow}
\usepackage{setspace}
\usepackage[title]{appendix}
\usepackage{mathrsfs} 
\usepackage{enumerate}

\usepackage{tikz}
\usepackage{tikz-cd}
\usepackage{wrapfig}
\usepackage{float}

\usetikzlibrary{matrix,arrows,decorations.pathmorphing}

\usepackage{amsthm}
\theoremstyle{plain}
\newtheorem{theorem}{Theorem}[section]

\theoremstyle{plain}
\newtheorem{lemma}[theorem]{Lemma}

\theoremstyle{plain}
\newtheorem{Prop}[theorem]{Proposition}

\theoremstyle{plain}
\newtheorem{corollary}[theorem]{Corollary}

\theoremstyle{definition}
\newtheorem{Def}{Definition}[section]
\theoremstyle{remark} 
\newtheorem{remark}{Remark}

\theoremstyle{definition}

\theoremstyle{definition}

\def\ed{\mathrm{d}}

\def\I{\mathcal I}

\def\J{\mathcal J}
\def\B{\mathcal B}
\def\K{\mathcal K}
\def\U{\mathcal U}

\def\F{\mathcal F}

\def\R{\mathbb R}

\def\P{\mathbb P}

\def\rank{{\rm rank}}

\def\W{\wedge}
\def\<{\langle}
\def\>{\rangle}

\def\GL{{\rm GL}}

\def\SL{{\rm SL}}

\def\lbb{[\![}
\def\rbb{]\!]}

\def\({\left(}
\def\){\right)}

\def\alg{{\rm alg}}

\def\bs{\boldsymbol}

\DeclareFontFamily{U}{MnSymbolC}{}
\DeclareSymbolFont{MnSyC}{U}{MnSymbolC}{m}{n}
\DeclareFontShape{U}{MnSymbolC}{m}{n}{
    <-6>  MnSymbolC5
   <6-7>  MnSymbolC6
   <7-8>  MnSymbolC7
   <8-9>  MnSymbolC8
   <9-10> MnSymbolC9
  <10-12> MnSymbolC10
  <12->   MnSymbolC12}{}
\DeclareMathSymbol{\lefthook}{\mathbin}{MnSyC}{'270}

\title[Rank 2 B\"acklund Transformations]{Rank 2 B\"acklund transformations of 
hyperbolic Monge-Amp\`ere systems}

\author{Yuhao Hu} 
\address{Jingdezhen, Jiangxi Province, P. R. China}
\email{YuhaoHu8@gmail.com}


\begin{document}

\maketitle

\begin{abstract}
		There are two main types of rank~2 B\"acklund transformations relating a pair of hyperbolic 
	Monge-Amp\`ere systems, which we call Type $\mathscr{A}$ and Type $\mathscr{B}$. 
		For Type $\mathscr{A}$, we completely determine a subclass whose local invariants satisfy a 
	specific but simple algebraic constraint; such B\"acklund transformations are parametrized by a finite 
	number of constants, whose cohomogeneity can be either 2, 3 or 4.
		In addition, we present an invariantly formulated condition that determines whether a generic 
	Type $\mathscr{B}$ B\"acklund transformation is one that, under suitable choices of local coordinates, 
	relates solutions of two PDEs of the form $z_{xy} = F(x,y,z,z_x,z_y)$ and preserves the $x,y$ variables on 
	solutions.
\end{abstract}

\section{Introduction}

	The study of B\"acklund transformations began in the 1880s when B\"acklund and Bianchi
discovered a relation between the so-called \emph{line congruences} and \emph{pseudospherical surfaces} 
in the three-dimensional Euclidean space. (For a modern treatment, see \cite{chern-terng}.)
	That particular discovery triggered a long search for pairs of partial differential equations (PDEs) whose solutions 
are related by a \emph{B\"acklund transformation}.
	The search is easily justified: having a B\"acklund transformation will allow one to find many solutions to a nonlinear
PDE by solving only ordinary differential equations (ODEs). 
	The study of B\"acklund transformations flourished in the mid-twentieth century, influencing differential 
geometry and mathematical physics. 
	However, it also gives the impression that pairs of PDEs that do admit a B\"acklund transformation seem to be isolated 
and rare compared with the classes of PDEs being considered. 
	The problem of determining the generality of B\"acklund transformations remains largely unsolved.

	Recently there seems to be a revival of interest in studying B\"acklund transformations. 
	A possible motivation for this is the emergence of a useful geometric formulation, which describes a B\"acklund 
transformation as a double fibration of smooth manifolds that carry certain geometric structures. 
	Using this formulation, one is essentially free from manipulating coordinate variables and can instead work with
quantities of intrinsic nature. 
 	This approach has already led to a number of new results, which include the following. 
	In \cite{C01}, using Cartan's methods, Clelland obtained a complete classification of homogeneous rank $1$ 
B\"acklund transformations relating a pair of hyperbolic Monge-Amp\`ere systems, where several analogues of 
the classical example of B\"acklund and Bianchi were found. 
	In \cite{AF12} and \cite{AF15}, Anderson and Fels developed a method of constructing B\"acklund transformations
using symmetry reduction, which also addressed a connection between B\"acklund transformations and Darboux
integrability, a direction previously explored by Clelland and Ivey in \cite{CI09}.
	Later, a local generality result appeared in \cite{HuBacklund1}, where the homogeneity assumption made 
in \cite{C01} was removed; the paper confirmed rarity of existence in a generic case and produced some new examples.

	In Section \ref{BasicDefNConstruction} of this article, we will recall the geometric formulation mentioned above.
	From the definition, a notion of `rank' naturally arises; it represents the fiber rank of the underlying double fibration, 
which also measures the abundance of solutions a B\"acklund transformation can produce from a given solution. 
	While many classical B\"acklund transformations have rank $1$, B\"acklund transformations of higher ranks deserve 
no less attention. 
	For instance, the classical Tzitz\'eica transformation, which was first discovered in Tzitz\'eica's study 
(\cite{Tz08}, \cite{Tz09}) of affine spheres and was more recently revisited by \cite{dunajski2002hyper} and \cite{wang2006tzitzeica}
from some new perspectives, 
is a $1$-parameter family of rank $2$ B\"acklund transformations. 
	In \cite{AF15}, it was pointed out that the pair of PDEs
	\[
		u_{xy} = \frac{\sqrt{1-{u_x}^2}\sqrt{1 - {u_y}^2}}{\sin u} \quad \text{and}\quad v_{xy} = 0,
	\]
which admits a rank $2$ B\"acklund transformation, doesn't admit any rank $1$ (real) B\"acklund transformation between 
their solutions.

	The objective of this article is to take a beginning step towards a systematic study of rank $2$ B\"acklund 
transformations relating a pair of hyperbolic Monge-Amp\`ere systems.
	Before delving into technicalities, one should be aware that some B\"acklund transformations 
of higher ranks may be constructed in obvious
ways from those of lower ranks (see Section \ref{HigherRankConstr}).
	Consequently, we will exclude from our study those rank $2$ B\"acklund transformations that can arise trivially from a 
$1$-parameter family of rank $1$ B\"acklund transformations.

	With these in mind, we address the equivalence problem, in the sense of \'Elie Cartan, for the B\"acklund 
transformations under consideration.
	At an early stage, we identify two main types of B\"acklund transformations, 
which we call Type $\mathscr{A}$ and Type $\mathscr{B}$ (Section \ref{MainTypes}).

	For generic Type $\mathscr{A}$ B\"acklund transformations, while their full generality remains to be determined,
we find that (Theorem \ref{caseA1rigidityThm}), if the relative invariants of the associated $G$-structure 
take value in a specific codimension 4 subspace (characterized by \eqref{caseA1ThmAssu}) inside a total space 
of 31 dimensions, then the resulting B\"acklund transformations are parametrized by a finite number of constants, 
whose structure equations are completely determined (see \eqref{caseA1oms} and \eqref{caseA1Hs}).
	Such B\"acklund transformations can only have cohomogeneity 2, 3 or 4, which is determined by the image of 
an intrinsically defined map $\Theta$ from the $7$-dimensional base manifold of a B\"acklund transformation to 
a $4$-dimensional space. 
	Particularly, in the cohomogeneity 2 case, which is the case when the image of $\Theta$ is a surface,
the B\"acklund transformation must be one that relates solutions of the PDE
\[
	(x+y)z_{xy} + 2\sqrt{z_xz_y} = 0,
\]
which appears at the beginning of Goursat's list of Darboux integrable equations,
with those of the PDE
\[
	(x+y) z_{xy} - ({W_0}(e^{z_x})+1)({W_{-1}}(-e^{z_y})+1) = 0,
\]
where $W_0, W_{-1}$ are the two real branches of the Lambert $W$ function.

	For generic Type $\mathscr{B}$ B\"acklund transformations, we obtain a criterion (Proposition \ref{caseB12DorbitProp})
for determining when, up to contact transformations, a B\"acklund transformation is one that relates a pair of PDEs 
of the form
		\[
			z_{xy} = F(x,y,z,z_x,z_y)\quad \text{and}\quad Z_{XY} = G(X,Y,Z,Z_X,Z_Y),
		\]
in such a way that $x = X$ and $y = Y$ on corresponding solutions.

Most calculations in this article are performed 
using Maple\texttrademark.

\section{B\"acklund Transformations}\label{BasicDefNConstruction}

	To geometrically define B\"acklund transformations, we need the notion of an \emph{exterior differential system}
and that of an \emph{integrable extension}.

\begin{Def}\label{EDSDef}
	An \emph{exterior differential system} (EDS) is a pair $(M,\I)$, where $M$ is a smooth manifold, 
	and $\I\subset \Omega^*(M)$ is an ideal that is closed under exterior differentiation.
\end{Def}

\begin{Def}\label{IntMfdDef}
	An \emph{integral manifold} of an EDS $(M,\I)$ is a submanifold 
	\[
		\iota: N\hookrightarrow M
	\]
	that satisfies $\iota^*\I = 0$. 
\end{Def}

The following correspondence is well-known (see \cite{BCG}):

	\begin{center}
		EDS $(M,\I; \alpha)$$\leftrightsquigarrow$ PDE system $\mathcal{E}$\\
		Integral manifolds $\leftrightsquigarrow$ Solutions
	\end{center}
where $\alpha\in \Omega^*(M)$ is an independence condition, which corresponds to the independent
variables on the PDE side.

\begin{Def}\label{IntExtDef}
	Let $(M,\I)$ be an EDS. A rank $k$ \emph{integrable extension} of $(M,\I)$ is a
	submersion $\pi: (N,\J)\rightarrow(M,\I)$ with fibre dimension $k$ such that, for each $p\in N$,
	there exists an open neighborhood $U$ of $p$ and a rank $k$ vector subbundle $\Theta\subset T^*U$ 
	that satisfy the following two conditions:
	\begin{enumerate}[(1)]
		\item{On $U$, $\J$ is algebraically generated by $\pi^*\I$ and the sections of $\Theta$;}
		\item{$\Theta^\perp\cap \ker \pi_*  = 0.$}
	\end{enumerate}
\end{Def}

	Let $\pi: (N,\J)\rightarrow (M,\I)$ be a rank $k$ integrable extension. Given any integral manifold 
$S\subset M$ of $(M,\I)$, it follows from the Frobenius theorem that $\pi^{-1}S$ is foliated by a $k$-parameter 
family of integral manifolds of $(N,\J)$. 
	Furthermore, restricted to any integral manifold of $(N,\J)$, $\pi$ is a local diffeomorphism onto its image, 
which is an integral manifold of $(M,\I)$.

\begin{Def}\label{BacklundDef}
	A \emph{B\"acklund transformation} $(N,\B;\pi,\bar\pi)$ relating two EDS $(M,\I)$ and $(\bar M, \bar \I)$ 
	is a double fibration as indicated in the diagram below, where $\pi$ and $\bar\pi$ are both integrable extensions. 
	In particular, when $(M,\I)$ and $(\bar M,\bar \I)$ are contact equivalent, $(N,\B;\pi,\bar\pi)$
	is called an \emph{auto-B\"acklund transformation} of $(M,\I)$.
			\begin{figure*}[h!]
				\begin{center}
					\begin{tikzcd}[column sep=tiny]
							    &   (N,\mathcal{B}) \arrow[ld,"\displaystyle\pi",swap]\arrow[rd,"\displaystyle{\bar{\pi}}"]  
							    			& 		\\
							(M,\I)  &  &(\bar M,\bar \I)
 					\end{tikzcd}
					\label{BacklundDiagram}
				\end{center}	
				\end{figure*}		 	
\end{Def}

	From the `Frobenius' argument above, it follows that a B\"acklund transformation allows one to 
obtain integral manifolds of $(\bar M,\bar\I)$ from a given integral manifold of $(M,\I)$, and vice versa,
for which only ODE techniques are needed.

\begin{Def}\label{BacklundRankDef}
	We say that a B\"acklund transformation $(N,\B;\pi,\bar\pi)$ has \emph{rank $r$} if both $\pi$ and $\bar\pi$, 
	as integrable extensions, have rank $r$.
\end{Def}

\subsection{Hyperbolic Monge-Amp\`ere systems}

	In terms of the above, this article is concerned with rank $2$ B\"acklund transformations. In addition, 
we will assume that the two EDS being B\"acklund-related are both hyperbolic Monge-Amp\`ere systems, 
defined as follows.

 \begin{Def}\label{MADef}
 	A \emph{hyperbolic Monge-Amp\`ere system} is an EDS $(M,\I)$, where $M$ is a $5$-dimensional 
	smooth manifold, and $\I$ is locally generated by
	a contact form $\theta$, its exterior derivative $\ed \theta$ and a two form $\Phi$ such that, at each point,
	all solutions $[\lambda,\mu]$ of the congruence
		\[
			(\lambda\ed\theta+\mu\Phi)^2\equiv 0 \mod \theta
		\]
	correspond to precisely two distinct points in  $\R\P^1$.	
 \end{Def}

	Given a hyperbolic Monge-Amp\`ere system $(M,\I)$, locally one can always choose a coframing 
$\bs\eta = (\theta,\eta^1,\ldots,\eta^4)$ such that 
	\[
		\I = \<\eta^0, \eta^1\W\eta^2,\eta^3\W\eta^4\>.
	\] 
(Here and in the following, $\<\cdots\>$ denotes the \emph{differential} ideal generated by 
the enclosed forms and their exterior derivatives; $\<\cdots\>_\alg$ will denote the \emph{algebraic} ideal
generated by the forms alone.)

	The unordered pair of rank $3$ vector bundles 
\[
	\Xi_{10}:= \lbb\theta,\eta^1,\eta^2\rbb,\qquad \Xi_{01}:= \lbb\theta,\eta^3,\eta^4\rbb
\]		
is intrinsically defined. They are called the \emph{characteristic systems} of $(M,\I)$. 
(We will use $\lbb\cdots\rbb$ to denote the vector bundle generated by the enclosed elements.)

	A particular class of hyperbolic Monge-Amp\`ere systems are those that correspond to 
second-order PDEs of the form
	\[
		z_{xy} = F(x,y,z, z_x, z_y).
	\]
	In EDS language, let $M\subset J^1(\R^2,\R)$ be an open subset on which the function $F = F(x,y,z,p,q)$ is defined; 
and let $\theta = \ed z - p \ed x  - q\ed y$ be the restriction to $M$ of the canonical contact form on $J^1(\R^2,\R)$; 
finally let 
 	\[
		\I = \<\theta,\ed \theta, \(\ed p  - F(x,y,z,p,q)\ed y\)\W\ed x\>_\alg.
	\] 	
	This class of hyperbolic Monge-Amp\`ere systems has an intrinsic characterization, as the following
proposition shows (see \cite{HuBacklund2} for a proof).

\begin{Prop}\label{genFGordonProp}
	Locally a hyperbolic Monge-Amp\`ere system corresponds to a PDE of the form
		\[
			z_{xy} = F(x,y,z,z_x,z_y)
		\]
	up to contact equivalence if and only if each of the characteristic systems $\Xi_{10}$ and $\Xi_{01}$ 
	admits a nontrivial first integral.
\end{Prop}

\subsection{B\"acklund Transformations of Higher Ranks}\label{HigherRankConstr}

	We now describe two constructions of B\"acklund transformations of ranks higher than $1$
from those of lower ranks. (c.f. \cite{CI05} and \cite{AF16}.)

\subsubsection{Construction from a family}\label{familyconstrsection}
	As a prototypical example, one can regard the classical $1$-parameter family of sine-Gordon transformations as 
a single rank $2$ B\"acklund transformation of the sine-Gordon equation; details are left to the reader.
	More generally, let $(M,\I)$ and $(\bar M,\bar \I)$ be two fixed EDS that are related by a smooth family of rank $1$ 
B\"acklund transformations $(N,\B_\lambda; \pi_\lambda,\bar\pi_\lambda)$ $(\lambda\in\R)$.
	By ``smooth" we mean that (a) there exist smooth submersions
	\[
		\pi: N\times \R \rightarrow M, \qquad \bar\pi: N\times \R\rightarrow \bar M,
	\]
such that 
	\begin{equation}\label{pilambdaetc}
		\pi_\lambda = \pi(\cdot, \lambda)  \quad \text{and}\quad \bar\pi_\lambda = \bar\pi(\cdot,\lambda),
	 \end{equation}
and (b) each point $p\in N$ has an open neighborhood $U$, on which there exist two families of nonvanishing 
$1$-forms $\theta_\lambda$ and $\bar\theta_\lambda$ $(\lambda\in \R)$, varying smoothly in $\lambda$, 
such that, for each fixed $\lambda$,
	\[
		\mathcal{B}_\lambda|_U = \<\pi_\lambda^*\I, \theta_\lambda\>_\alg = 
		\<{\bar\pi_\lambda}^*\bar\I, \bar\theta_\lambda\>_\alg.
	\]

	Consider the space $\hat N:= N\times \R$, and let  $\lambda$ be the coordinate on the $\R$-component. 
Let  $\rho$ and $\sigma$ be the obvious projections shown in 
the diagram below.
	\begin{figure*}[h!]
				\begin{center}
					\begin{tikzcd}[column sep=small]
							    &  \hat N\arrow[ld,"\displaystyle\rho",swap]\arrow[rd,"\displaystyle{{\sigma}}"]  
							    					& 		\\
							N  &  &\R
 					\end{tikzcd}
				\end{center}	
	\end{figure*}

	Now define an algebraic ideal $\hat\B$ in $\Omega^*(\hat N)$ by
	\[
		\hat \B_{(p,\mu)}:= \<\rho^*\B_{\mu}, \ed \lambda\>_\alg,\quad \forall (p,\mu)\in \hat N.
	\]
	This $\hat \B$ is essentially generated by $\ed \lambda$  and the algebraic generators of $\B_\lambda$ 
in which $\lambda$ shall be viewed as a variable rather than a constant.

\begin{Prop}
	$(\hat N,\hat \B;\pi,\bar \pi)$ is a rank $2$ B\"acklund transformation relating $(M,\I)$ and $(\bar M,\bar \I)$.
\end{Prop}

\emph{Proof.} For any $(p,\mu)\in \hat N$, let $U\subset N$ and $\theta_\mu\in \Omega^1(U)$ be as above. 
By assumption, we have
	\[
		\B_\mu = \<\pi_\mu^*\I, \theta_\mu\>_\alg
	\]
and
	\[
		\hat \B_{(p,\mu)} = \<\rho^*\pi_\mu^*\I, \rho^*\theta_\mu, \ed \lambda\>_\alg
					 = \<\pi^*\I, \rho^*\theta_\mu, \ed\lambda\>_\alg.
	\]
	The last equality holds because for any $1$-form $\omega\in \Omega^1(M)$, the difference 
$\pi^*\omega  - \rho^*\pi_\mu^*\omega$ evaluated at $(p,\mu)\in \hat N$ is a scalar multiple of $\ed \lambda$.
	Furthermore, as $(p,\mu)$ varies, $\rho^*\theta_\mu$ gives rise to a $1$-form $\theta$ defined 
on $U\times\R\subset \hat N$.
	Thus, on $U\times\R$ we have
		\[
			\hat \B = \<\pi^*\I, \theta,\ed \lambda\>_\alg,
		\]
which is differentially closed, because $\ed\theta$ splits into two parts, one belonging to $\pi^*\I$ and the other 
being a multiple of $\ed\lambda$.
	
	It follows that $\pi: (\hat N,\hat\B)\rightarrow (M,\I)$ is an integrable extension. 
The case for $\bar\pi$ is similar.
\qed
\vskip 2mm

	In the same way, one can obtain a B\"acklund transformation of rank $k+r$ from a $k$-parameter family of rank $r$ 
B\"acklund transformations relating a fixed pair of EDS.

\subsubsection{Construction by composition}

	Consider two B\"acklund transformations, one relating $(M_1, \I_1)$ and $(M_2, \I_2)$, the other relating 
$(M_2, \I_2)$ and $(M_3, \I_3)$, as the following diagram shows.
	\begin{figure*}[h!]
		\begin{center}
			\begin{tikzcd}[column sep=tiny]
				&   (N_1,\mathcal{B}_1) \arrow[ld,"\displaystyle\pi_1",swap]\arrow[rd,"\displaystyle{{\pi}_2}"]  & 	
					&   (N_2,\mathcal{B}_2) \arrow[ld,"\displaystyle\pi_3",swap]\arrow[rd,"\displaystyle{{\pi}_4}"]  	
				\\
				(M_1,\I_1)  &  &(M_2,{\I}_2) &&(M_3, \I_3)
 			\end{tikzcd}
		\label{Backlund}
		\end{center}	
	\end{figure*}

	The Whitney sum of the fiber bundles $\pi_2: N_1\rightarrow M_2$ and $\pi_3: N_2\rightarrow M_2$, 
denoted as $N_1\oplus N_2$, admits two submersions
	\[
		p_1: N_1\oplus N_2\rightarrow N_1, \qquad p_2: N_1\oplus N_2\rightarrow N_2,
	\] 
satisfying $\pi_2\circ p_1 = \pi_3\circ p_2$.	
Let $\B$ denote the differential ideal on $N_1\oplus N_2$ algebraically generated by $p_1^*\B_1$ and $p_2^*\B_2$.

\begin{Prop}\label{composingBacklund}
	$(N_1\oplus N_2, \B; \pi_1\circ p_1, \pi_4\circ p_2)$ is a B\"acklund transformation relating $(M_1, \I_1)$ 
	and $(M_3, \I_3)$.
\end{Prop} 

\emph{Proof}. 
Considering Lemma \ref{intextcomp} below, it suffices to show that $p_1$ and $p_2$ are integrable extensions.  
By the assumption, there exist $1$-forms $\alpha_i,\beta_j$ such that 
	\[
		\B_1= \<\alpha_1,\ldots,\alpha_k, \pi_2^*\I_2\>_\alg,\qquad 
	\B_2 =\<\beta_1,\ldots,\beta_\ell, \pi_3^*\I_2\>_\alg.
	\]
	 Thus, by construction,
	 \[
	 	\B = \< p_1^*\alpha_1,\ldots,p_1^*\alpha_k, p_2^*\beta_1,\ldots,p_2^*\beta_\ell, 
			p_1^*\pi_2^*\I_2\>_\alg.
	\]		
	Note that $p_1^*\pi_2^*\I_2$ is the same as
					$p_2^*\pi_3^*\I_2$.  
					
	Now suppose that ${\bs v}\in T(N_1\oplus N_2)$ is tangent to a fiber of $p_1$. 
	It follows that ${p_2}_*\bs v$ is tangent to a fiber of $\pi_3$.
	If $p_2^*(\beta_j)({\bs v}) = \bs 0$ for all $j= 1,\ldots,\ell$, it is necessary that ${p_2}_*{\bs v} = \bs 0$, 
because $\pi_3$ is an integrable extension. 
	Since $p_2$, restricted to each fiber of $p_1$, is an immersion, $\bs v$ must vanish. 
	This proves that $p_1$ is an integrable extension. The case for $p_2$ is similar. \qed

\begin{lemma}\label{intextcomp}
	The composition of two integrable extensions is an integrable extension. 
\end{lemma}

\emph{Proof}. 
Suppose that
	\[
		\pi_1:(M, \I)\rightarrow (N, \J), \quad \pi_2:(N, \J)\rightarrow (P,\K)
	\]
are integrable extensions of ranks $q$ and $p$, respectively.
	We prove that
	\[
		\pi_2\circ\pi_1: (M,\I)\rightarrow (P,\K)
	\]
is an integrable extension of rank $p+q$.
   
By definition, locally $\J$ is algebraically generated by $\pi_2^*\K$ and some $1$-forms $\alpha_1,\ldots,\alpha_p$. 
Thus, there exist $q$ $1$-forms $\beta_1,\ldots,\beta_q$ such that $\I$ is algebraically generated by 
	\[ 
		(\pi_2\circ\pi_1)^*\K, \pi_1^*\alpha_1,\ldots,\pi_1^*\alpha_p, \beta_1,\ldots,\beta_q.
	\]
For $\pi_2\circ\pi_1$, the first condition in Definition \ref{IntExtDef} is clearly satisfied.
To verify the second condition, suppose that there exist constants $c_i, f_j$ such that 
	\begin{equation}\label{compIntExtEq}
		\left(\sum_{i = 1}^p c_i \pi_1^*\alpha_i + \sum_{j = 1}^q f_j \beta_j\right)({\bs v}) = 0
	\end{equation}	
for any ${\bs v}\in T_xM$ satisfying ${\pi_2}_*({\pi_1}_*(\bs v)) = 0$.  
Since $\pi_1$ is an integrable extension, each $T_xM$ is a direct sum:
	\[	
		T_xM = V_1\oplus V_2 :=\ker_x({\pi_1}_*) \oplus \{\beta_1,\ldots,\beta_q\}_x^\perp.
	\]	 
In order for \eqref{compIntExtEq} to hold on $V_1$, all $f_j$ must vanish. In order for \eqref{compIntExtEq} to hold on 
$V_2\cap \ker((\pi_2\circ\pi_1)_*)$, $c_i$ must all vanish, because the restriction of ${\pi_1}_*$ to $V_2$ is 
a linear isomorphism. 
This completes the proof.
\qed

\begin{corollary}
	Being B\"acklund-related is an equivalence relation for exterior differential systems. 
\end{corollary}	

\emph{Proof.} 
Symmetry and reflexivity are obvious from the definition. 
Transitivity is just Proposition \ref{composingBacklund}.
\qed

\section{Genericity Assumptions and Structure Reduction}\label{MainTypes}

	From now on, $(N,\B; \pi,\bar \pi)$ will denote a rank $2$ B\"acklund transformation 
relating a pair of hyperbolic Monge-Amp\`ere systems $(M,\I)$ and $(\bar M,\bar\I)$. 
	As a hyperbolic EDS, $(N,\B)$ is of class $s=3$ in the sense of \cite{BGH1}; it has an unordered 
pair of characteristic bundles $\chi_{10}$ and $\chi_{01}$. 
	However, in our analysis below, we will regard $\chi_{10}$ and $\chi_{01}$ as an \emph{ordered} pair,
for better clarity.

	Since $(M,\I)$ is hyperbolic Monge-Amp\`ere, for every point $x\in M$, there exists (according to \cite{BGG}) 
an open neighborhood $\U\subset M$ of $x$ and an \emph{adapted} coframing
$(\theta,\eta^1,\ldots,\eta^4)\in \F^*(\U)$, which by definition satisfies
	\[		
		\theta\W(\ed\theta)^2\ne 0,\quad 
		\I = \<\theta, \eta^1\W\eta^2, \eta^3\W\eta^4\>_\alg.
	\]						
Similarly, for every point $\bar x\in \bar M$, there exists a neighborhood $\bar\U$ of $\bar x$ and 
an adapted coframing $(\bar\theta,\bar\eta^1,\ldots,\bar\eta^4)\in \F^*(\bar\U)$.

\begin{Def}\label{transversalDef}
	Let $\pi:N\rightarrow M$ be a submersion of smooth manifolds. A vector subbundle $J\subset T^*N$ is said to be
 \emph{transversal} to $\pi$ if, at each point $p\in N$, 
 	\[
		(J_p)^\perp \cap \ker_p(\pi_*) = {\bf 0}.
	\]
\end{Def}
 \vskip 2mm

	By Definitions \ref{BacklundDef} and \ref{transversalDef},
there exists a rank $2$ subbundle $J\subset B^1$ 
that is transversal to $\pi$, which satisfies
					\[
						B^1  = J\oplus \pi^*(I^1).
					\]				
(Here and in the following, we use $I^k$ to denote the vector bundle
generated by the degree $k$ part of a differential ideal $\I$.)
Similarly, there exists a rank $2$ subbundle $\bar J\subset B^1$, transversal to $\bar\pi$, satisfying
					\[
						B^1 = \bar J\oplus\bar \pi^*(\bar I^1).
					\]
In particular, $B^1\subset T^*N$ has rank $3$.

	Next, we will set up the \emph{equivalence problem}, in the sense of \'Elie Cartan (see \cite{GardnerEquiv}),
for rank $2$ B\"acklund transformations relating two hyperbolic Monge-Amp\`ere systems. 
As we proceed, we will state two \emph{genericity conditions} to help identify the cases of interest.

We begin by stating the following.
\begin{quote}
	\emph{First genericity condition:}  \[\pi^*I^1\cap\bar\pi^*\bar I^1 = \bs 0.\] 
\end{quote}

	Assuming this condition, for each $p\in N$, there exists an open neighborhood $U$ of $p$, 1-forms 
$\gamma\in \Omega^1(U)$, $\theta\in\I^1|_{\pi(U)}$ and $\bar\theta\in\bar \I^1|_{\bar\pi(U)}$
such that $\pi^*\theta, \bar\pi^*\bar\theta, \gamma$ span $B^1|_U$. (Here $\I^k:= \I\cap \Omega^k(M)$, etc.)
	Moreover, $\gamma$ can be chosen in such a way that 
$\pi^*\bar\theta$ and $\gamma$  (respectively, $\pi^*\theta$ and $\gamma$) are linearly independent when 
pulled back to each fiber of 
$\pi$ (respectively, $\bar\pi$).

	By shrinking $U$, if needed, we can extend $\theta$ to a coframing $(\theta,\eta^1,\ldots,\eta^4)$
on $\pi(U)$ that is adapted to the Monge-Amp\`ere ideal $\I$; similarly extend $\bar\theta$ to a coframing
$(\bar\theta,\bar\eta^1,\ldots,\bar\eta^4)$ on $\bar\pi(U)$, adapted to $\bar \I$. 
Doing this, we obtain a coframing on $U\subset N$:
	\[
		(\pi^*\theta, \bar\pi^*\bar\theta, \gamma, \pi^*\eta^1, 
			\pi^*\eta^2, \pi^*\eta^3, \pi^*\eta^4).
	\] 
Dropping the pullback symbols for clarity, we have that
	\begin{equation*}
		\B|_U = \<\theta, \bar\theta, \gamma, \eta^1\W\eta^2, \eta^3\W\eta^4\>_{\rm alg}
			   = \<\theta,\bar\theta,\gamma,\bar\eta^1\W\bar\eta^2, \bar\eta^3\W\bar\eta^4\>_{\rm alg}.
	\end{equation*}

Now we state the following.
 	\begin{quote}
 	 	\emph{Second genericity condition:} 
		\[
			\lambda \ed\theta + \mu\ed\bar\theta \equiv 0 \mod \theta,\bar\theta, \gamma
		\]
		if and only if $\lambda = \mu = 0$.
		 (Note that this condition does not depend on the choice of $\theta,\bar\theta$ and $\gamma$.)
	\end{quote}

Assuming both genericity conditions, we have that
			\[
				\lbb \ed\theta,\ed\bar\theta\rbb 
					\equiv
				\lbb\eta^1\W\eta^2, \eta^3\W\eta^4\rbb
					\equiv 
				\lbb\bar\eta^1\W\bar\eta^2, \bar\eta^3\W\bar\eta^4\rbb 
				\quad \mod\theta,\bar\theta,\gamma.
			\]
for the three rank $2$ vector bundles involved.					
In addition, 
			\[
				\B = \<\theta,\bar\theta,\gamma,\ed\theta,\ed\bar\theta\>_{\rm alg}.
			\]

The choice of a coframing on $U$ above has some ambiguity in it.
To refine the choice, we start with an adapted coframing $(\theta,\eta^1,\ldots,\eta^4)$ on $\pi(U)\subset M$ 
that satisfies the extra condition
			\[
				\ed\theta\equiv \eta^1\W\eta^2+\eta^3\W\eta^4 \quad \mod\theta.
			\]
The same congruence then holds for the pullbacks of $\theta,\eta^1,\ldots,\eta^4$ by $\pi$.

Once we make the change of notations 
\[
	\pi^*\theta\mapsto \omega^0, \quad \bar\pi^*\bar\theta\mapsto \bar\omega^0, 
	\quad \pi^*\eta^i\mapsto \omega^i,
\]
we obtain the following congruences on $U\subset N$:
				\begin{align}
					\ed\omega^0
						&\equiv \omega^1\W\omega^2+\omega^3\W\omega^4 
						&& \mod \omega^0,	\label{hg2omega0}\\
					\ed\bar\omega^0
						&\equiv A_1\omega^1\W\omega^2+A_2\omega^3\W\omega^4  
						&&\mod \omega^0,\bar\omega^0, \gamma,	\label{hg2bomega0}\\
					\ed\gamma
						&\equiv A_3\omega^1\W\omega^2+A_4\omega^3\W\omega^4
						&&\mod\omega^0,\bar\omega^0,\gamma,	\label{hg2gamma}
				\end{align}
where $A_1,\ldots,A_4$ are functions defined on $U$.
Here $A_1,A_2$ are nonvanishing because $\bar\omega^0\W(\ed\bar\omega^0)^2\ne 0$. 
Moreover, $A_1\ne A_2$, by the second genericity assumption.

	Next, we perform the following steps successively:
\begin{enumerate}[\it Step 1.]
	\item{add multiples of $\omega^0,\bar\omega^0$ to $\gamma$ to arrange that $A_3 = A_4 = 0$.}
	\item{add multiples of $\omega^0$ to $\omega^1,\ldots,\omega^4$ such that
		\eqref{hg2bomega0} still holds when reduced modulo only $\bar\omega^0$ and $\gamma$. }
	\item{scale $\bar \omega^0$ to arrange that $A_1 = 1$.}
	\item{depending on the sign of $A_2$, scale $\omega^3$ and $\omega^0$
		 to put \eqref{hg2omega0} and \eqref{hg2bomega0} in the form
 		\begin{align*}
			\ed\omega^0
				&\equiv A\omega^1\W\omega^2 +\phantom{\epsilon} 
					\omega^3\W\omega^4 \quad\mod \omega^0,
				\\
			\ed\bar\omega^0 
				&\equiv \omega^1\W\omega^2 +\epsilon A \omega^3\W\omega^4 
					\quad \mod \bar\omega^0, \gamma,
		\end{align*}
		where $\epsilon = \pm 1$, $A>0$, and  $A^2\ne \epsilon$. }
	\item{if needed, swap $(\omega^1,\omega^2)$ with $(\omega^3,\omega^4)$, and
		then scale $\omega^0$ and $\bar\omega^0$ so that the characteristic bundles of $(N,\B)$ are
		\begin{equation}\label{orderedchar}
			\chi_{10} = \lbb\omega^0,\bar\omega^0,\gamma,\omega^1,\omega^2\rbb,
			\qquad
			\chi_{01} = \lbb\omega^0,\bar\omega^0,\gamma,\omega^3,\omega^4\rbb.
		\end{equation}
		}
	\item{add multiples of $\gamma$ to $\omega^3,\omega^4$ such that
		\begin{align}
			\ed\omega^0
				&\equiv A\omega^1\W\omega^2 + \phantom{\epsilon} \omega^3\W\omega^4 
					+ (B_3\omega^3+B_4\omega^4)\W\gamma \quad\mod \omega^0,	
					\label{generalomega0}
				\\
			\ed\bar\omega^0 
				&\equiv \omega^1\W\omega^2 + \epsilon A \omega^3\W\omega^4 
					+(B_1\omega^1+B_2\omega^2)\W\gamma \quad \mod \bar\omega^0,
					\label{generalbomega0}
			\end{align}
		for some functions $B_1,\ldots,B_4$ defined on $U$. Note that there cannot be an 
		$(\omega^0\W\gamma)$-term in $\ed\bar\omega^0$, because 
		$\bar \omega^0\W(\ed\bar\omega^0)^3 =0$.}
\end{enumerate}

Finally, there exist functions $C_0, C_i, D_i$ $(i = 1,\ldots, 4)$ on $U$ such that
		\begin{equation}
			\ed\gamma \equiv C_0\omega^0\W\bar\omega^0 +C_i\omega^i \W\omega^0 
						+D_i \omega^i\W\bar\omega^0 \quad \mod \gamma. \label{generaldgamma}
		\end{equation}

\begin{remark}
	The parameter $\epsilon$ has a geometric meaning:
	it indicates whether $(\ed\omega^0)^2$ and $(\ed\bar\omega^0)^2$ determine the same or the 
	opposite orientation(s) on the rank $4$ distribution $(B^1)^{\perp}$ on $U$.
\end{remark}

\begin{Def}
	We say that a rank $2$ B\"acklund transformation relating a pair of hyperbolic Monge-Amp\`ere systems is of 
	\begin{itemize}
	\item{\emph{Type $\mathscr{A}$}, if both genericity conditions hold;}
	\item{\emph{Type $\mathscr{B}$}, if only the first genericity condition holds.}
	\end{itemize}
\end{Def}

\section{Type $\mathscr{A}$}\label{typeAsection}

\begin{Def}
	 Let $(N,\B;\pi,\bar \pi)$ be a Type $\mathscr{A}$ rank $2$ B\"acklund transformation
	 with an ordered pair of characteristic systems $\chi_{10}$, $\chi_{01}$ for $(N,\B)$.
	A coframing
	 $(\omega^0,\bar\omega^0,\gamma,\omega^1,\ldots,\omega^4)$ 
	 defined on an open subset $U\subset N$ is said to be \emph{$0$-adapted}
	 if it satisfies 
 	\begin{equation}
 		\lbb\omega^0\rbb = \pi^*I^1, \quad 
		\lbb\bar\omega^0\rbb = \bar\pi^*\bar I^1, \quad 
		\lbb\omega^0,\bar\omega^0,\gamma\rbb = B^1,	\label{rank2zeroAdBundles}
	\end{equation}
	the condition \eqref{orderedchar} and 
	the congruences \eqref{generalomega0}, \eqref{generalbomega0}, \eqref{generaldgamma} 
	for a parameter $\epsilon  = \pm 1$ and some functions $A$ ($A > 0$ and $A^2\ne \epsilon$), 
	$B_1,\ldots, B_4$, $C_0, \ldots, C_4$, $D_1,\ldots, D_4$ defined on $U$. 
\end{Def}

\begin{lemma}\label{genAssu123Gstr}	
	Given a Type $\mathscr{A}$ rank $2$ B\"acklund transformation $(N,\B;\pi,\bar\pi)$, locally
	the $0$-adapted coframings are precisely the local sections of a $G$-structure $\mathcal{G}$ on $N$, 
	where $G\subset \GL(7,\R)$ is the subgroup generated by matrices of the form
			\begin{equation}					
						\begin{alignedat}{1}
							g  =& 	\left(
									\begin{array}{ccccc}
										\det(\bs b) &0&0&0&0\\	
										0&\det(\bs a)&0&0&0\\	
										0&0&c&0&0\\
										0&0&0&\bs a&0\\	
										0&0&0&0&\bs b\\	
									\end{array}
								\right), 
								\\
							& \det(\bs a) = \det(\bs b)\ne 0, \quad c\ne 0\in \R,\\
							&\quad \bs a = (a_{ij}),\bs b = (b_{ij})\in \GL(2,\R).
						\end{alignedat}						 
				\label{trans}
			\end{equation}
\end{lemma}

\emph{Proof}. Given a $0$-adapted coframing $\bs\omega$ defined on $U\subset N$, 
		it is easy to check that $\bs\omega\cdot g := g^{-1}\bs\omega$ ($\bs\omega$ viewed as a column) remains
		a $0$-adapted coframing for any $g: U\rightarrow G$.
		
		Conversely, by \eqref{orderedchar}, changing from one 
		$0$-adapted coframing to another must preserve the order of the characteristic bundles
		$\lbb\omega^0,\bar\omega^0,\gamma,\omega^1,\omega^2\rbb$
		and $\lbb\omega^0,\bar\omega^0,\gamma,\omega^3,\omega^4\rbb$.
		Thus, by \eqref{rank2zeroAdBundles}, 
		if $\bs\omega$ and $\hat{\bs\omega}$, both defined on $U\subset N$, are 
		two $0$-adapted coframings, then there exists a function $g:U\rightarrow \GL(7,\R)$ of the form
		\[
			g = \left(
									\begin{array}{ccc}
										 \Psi &0&0\\	
										*&\bs a&0\\	
										*&0&\bs b\\	
									\end{array}
								\right)\quad (\bs a,\bs b: U\rightarrow \GL(2,\R)),
		\]
		where $\Psi: U\rightarrow\GL(3,\R)$ is lower triangular with vanishing $(2,1)$-entry, 
		such that $\hat{\bs\omega} = \bs\omega\cdot g$.
		In order for $\bs\omega$ and $\hat{\bs\omega}$ to both satisfy 
		\eqref{generalomega0}, \eqref{generalbomega0} and \eqref{generaldgamma},
		$g$ must be of the form \eqref{trans}.		\qed

\begin{lemma}\label{GenAss123CharSysLemma}		
	Let $\bs\omega = (\omega^0,\bar\omega^0,\gamma,\omega^1,\ldots,\omega^4)$ be any
	$0$-adapted coframing of a Type $\mathscr{A}$ B\"acklund transformation $(N,\B;\pi,\bar\pi)$, defined
	on $U\subset N$.
	\begin{enumerate}[\rm (i)]
		\item{Let $\Xi_{10}=\lbb\theta,\eta^1,\eta^2\rbb$ and $\Xi_{01}=\lbb\theta,\eta^3,\eta^4\rbb$
 			be the restrictions to $\pi(U)$ of the characteristic systems of $(M,\I)$. 
			Their pullbacks via $\pi$ are, up to ordering, 
				\[\lbb\omega^0, \omega^1,\omega^2\rbb, \qquad 
				\lbb\omega^0, \omega^3-B_4\gamma, \omega^4+B_3\gamma \rbb.\]}
		\item{Let $\bar \Xi_{10}= \lbb\bar\theta,\bar\eta^1,\bar\eta^2\rbb$ and 
			$\bar\Xi_{01} = \lbb\bar\theta,\bar\eta^3,\bar\eta^4\rbb$ be the restrictions to $\bar\pi(U)$
			of the characteristic systems of $(\bar M,\bar \I)$. 
			Their pullbacks via $\bar\pi$ are, up to ordering, 
			\[
				\lbb\bar\omega^0,\omega^3,\omega^4\rbb,\qquad \lbb\bar\omega^0, 
				\omega^1-B_2\gamma, \omega^2+B_1\gamma\rbb.
			\]	}
\end{enumerate}			
\end{lemma}

\emph{Proof}. We only prove (i); the proof of (ii) is similar. 
		
		By \eqref{generalomega0}, the \emph{retracting space} (see \cite{BCG})
		of $\omega^0$ is
		\[
			\mathcal{C}(\<\omega^0\>) = \lbb\omega^0, \omega^1, \omega^2, 
		\omega^3 - B_4\gamma, \omega^4+B_3\gamma\rbb.
		\] 
		Because 
		\[
			\pi^*(\I)\subset \B, \quad  \pi^*(\lbb\theta\rbb) = \lbb\omega^0\rbb, \quad
			\mathcal{C}(\I) = \mathcal{C}(\<\theta\>),
		\]	 
		the ideal generated by
		$\pi^*(\I)$ must be contained in the \emph{intersection} of $\B$ and the algebraic ideal generated by 
		the sections of $\mathcal{C}(\<\omega^0\>)$, 
		which is 
		\begin{equation}
			\<\omega^0, \omega^1\W\omega^2, 
						(\omega^3-B_4\gamma)\W(\omega^4+B_3\gamma)\>_{\rm alg}. 
			\label{rank2InvCharsys}
		\end{equation}
		It follows that the ideal generated by
		$\pi^*(\I)$ is \emph{equal} to $\B\cap \mathcal{C}(\<\omega^0\>)$. 
		Comparing the characteristic systems
		leads to (i). \qed
\\

	By Lemma \ref{GenAss123CharSysLemma}, the $1$-forms 
$\omega^0, \omega^1,\omega^2,\omega^3-B_4\gamma, \omega^4+B_3\gamma$ are 
$\pi$-semibasic, meaning that the pullback of each of these forms to each fiber of $\pi$ is zero. 
	
	Moreover, the EDS \eqref{rank2InvCharsys} is invariant along the fiber directions of $\pi$. 
To be more precise, let 
	\[
		X_0, \bar X_0, X_\gamma, X_1,X_2,X_3,X_4
	\]
be the vector fields on $U\subset N$ that are dual to the $1$-forms  
	\[
		\omega^0,\bar\omega^0,\gamma,\omega^1,\omega^2, 
				\omega^3-B_4\gamma,\omega^4+B_3\gamma;
	\] 
and let 
	\[
		Y_0, \bar Y_0, Y_\gamma, Y_1,Y_2,Y_3,Y_4
	\]
be the vector fields dual to 
	\[
		\omega^0,\bar\omega^0,\gamma,\omega^1-B_2\gamma, 
				\omega^2+B_1\gamma,\omega^3,\omega^4.
	\] 
Lemma \ref{GenAss123CharSysLemma} implies that any $0$-adapted coframing defined on $U\subset N$ 
must satisfy the `invariance property' below, where $\mathcal{L}$ stands for the Lie derivative, and $\Gamma(U,K)$ stands for the
space of sections of a vector bundle $K$ over $U$.
	\begin{quote}
		\emph{Invariance property:} 
		
		When $K$ is either $\lbb\omega^0,\omega^1,\omega^2\rbb$ or
			$\lbb\omega^0,\omega^3-B_4\gamma,\omega^4+B_3\gamma\rbb$,
			\begin{equation}
				\sigma\in \Gamma(U,K)
				\Rightarrow
				\mathcal{L}_{\bar X_0}\sigma, ~\mathcal{L}_{X_\gamma}\sigma\in \Gamma(U,K);
				\label{InvProp1}
			\end{equation}
			
		when $K$ is either $\lbb\bar\omega^0,\omega^3,\omega^4\rbb$ or 
			$\lbb\bar\omega^0, \omega^1-B_2\gamma, \omega^2+B_1\gamma\rbb$,
			\begin{equation}
				\sigma\in \Gamma(U,K)
				\Rightarrow
				\mathcal{L}_{Y_0} \sigma, ~ \mathcal{L}_{Y_\gamma} \sigma\in \Gamma(U,K).
					\label{InvProp2}
			\end{equation}
	\end{quote}		
Conversely, we have the following proposition.

\begin{Prop}\label{allin7mfd}
	Let $U$ be a $7$-manifold. Suppose that there exists a coframing 
	$\bs\omega = (\omega^0,\bar\omega^0,\gamma,\omega^1,\ldots,\omega^4)$, 
	a constant $\epsilon = \pm 1$ and functions $A, B_i, C_0, C_i, D_i$ 
	$(i = 1,\ldots, 4; \epsilon = \pm 1; A>0, A^2\ne \epsilon)$
	defined on $U$ satisfying \eqref{generalomega0}, \eqref{generalbomega0} and \eqref{generaldgamma}. 
	If, in addition, $\bs\omega$ satisfies the invariance property described by \eqref{InvProp1} 
	and \eqref{InvProp2}, then $\bs\omega$ is $0$-adapted to a rank $2$ B\"acklund transformation 
	of Type $\mathscr{A}$.
\end{Prop}

\emph{Proof}. The proof is in three steps.

First, let
\begin{align*} 
\K:=&\<\omega^0,\omega^1,\omega^2,\omega^3-B_4\gamma,\omega^4+B_3\gamma\>,\\
\bar\K:=&\<\bar\omega^0,\omega^1 - B_2\gamma,\omega^2+B_1\gamma, \omega^3,\omega^4\>
\end{align*}
be two EDS on $U$ obtained from the coframing $\bs\omega$, where $B_i$ are the coefficients
occurring in \eqref{generalomega0} and \eqref{generalbomega0}. We claim that 
both $\K$ and $\bar \K$ are Frobenius. Assuming this, we can construct
 $\pi:U\rightarrow M$ as the quotient
of $U$ by the $2$-dimensional
leaves of $\K$, and similarly obtain $\bar\pi: U\rightarrow \bar M$. 

To see that $\K$ is Frobenius, note that the vector fields $\bar X_0$ and $X_\gamma$ annihilate
all degree $1$ generators of $\K$; thus, it suffices to show that their Lie bracket 
$[\bar X_0, X_\gamma]$ satisfies the same
property. We verify this for $\omega^1$. In fact,
	\[
		\omega^1([\bar X_0, X_\gamma]) = -\ed\omega^1(\bar X_0, X_\gamma) 
		= -(\bar X_0\lefthook \ed\omega^1)(X_\gamma).
	\]
By \eqref{InvProp1}, $\bar X_0\lefthook \ed\omega^1 = \mathcal{L}_{\bar X_0}\omega^1$ is a linear combination of $\omega^0,\omega^1,\omega^2$; hence, 
	\[
		(\bar X_0\lefthook \ed\omega^1)(X_\gamma) = 0.
	\]
The other cases, including the case for $\bar\K$, are similar. 

Second, we need to show that the systems
	\begin{align*} 
		\J:=& \<\omega^0,\omega^1\W\omega^2, (\omega^3 - B_4\gamma)\W(\omega^4+B_3\gamma)\>,\\
		\bar\J:=& \<\bar\omega^0, (\omega^1 - B_2\gamma)\W(\omega^2+B_1\gamma),
		\omega^3\W\omega^4\>
	\end{align*}
descend to $M$ and $\bar M$, respectively. To see this, note that \eqref{generalomega0}
implies that 
$\mathcal{L}_{\bar X_0} \omega^0, \mathcal{L}_{X_\gamma}\omega^0$ are both
multiples of $\omega^0$; hence, there exists a contact form 
$\eta^0\in \Omega^1(M)$ satisfying 
$\lbb\pi^*\eta^0\rbb = \lbb\omega^0\rbb$. By the invariance property, there exists
a coframing $\bs\eta = (\eta^0,\ldots,\eta^4)$ on $M$ such that 
\begin{align*}
	\pi^*\lbb\eta^0,\eta^1,\eta^2\rbb &=  \lbb\omega^0,\omega^1,\omega^2\rbb,\\
	\pi^*\lbb\eta^0,\eta^3,\eta^4\rbb &= \lbb\omega^0, \omega^3 - B_4\gamma, \omega^4+B_3\gamma\rbb.
\end{align*}
 It follows that $\I = \<\eta^0, \eta^1\W\eta^2,\eta^3\W\eta^4\>_\alg$ is hyperbolic Monge-Amp\`ere
 and that $\J$ is generated by $\pi^*\I$. The case for $\bar \J$ is similar. 

Third, consider the differential ideal $\B = \<\omega^0, \bar\omega^0, \gamma,\omega^1\W\omega^2,\omega^3\W\omega^4\>_\alg$ on $U$. 
It is easy to verify by definition that $(U,\B;\pi,\bar\pi)$ is a 
Type $\mathscr{A}$ rank $2$ B\"acklund transformation relating $(M,\I)$ and $(\bar M,\bar\I)$.
Moreover, the coframing $\bs\omega$ is $0$-adapted to this B\"acklund transformation 
with the ordered pair of characteristic systems $\chi_{10} = \lbb\omega^0,\bar\omega^0,\gamma,
\omega^1,\omega^2\rbb$ and $\chi_{01}=\lbb\omega^0,\bar\omega^0,\gamma,
\omega^3,\omega^4\rbb$.
 \qed
\\

Proposition \ref{allin7mfd} implies that the geometry of Type $\mathscr{A}$ B\"acklund transformations
is completely contained in a $G$-structure on a $7$-manifold.

     Let $\mathcal{G}$ be as in Lemma \ref{genAssu123Gstr}.
     By a simple calculation, one can show that, under the transformation 
      $u\mapsto u\cdot g = g^{-1}u\in \mathcal{G}$, where
      $g$ is as in $\eqref{trans}$, the coefficients $B_1,\ldots, B_4$ in 
      \eqref{generalomega0} and \eqref{generalbomega0} transform by
      					\begin{equation*}
						\begin{alignedat}{1}
						\left(
							\begin{array}{c}
								B_1\\
								B_2
							\end{array}
						\right)(u\cdot g)& = 
								\frac{c}{\det(\bs a)}\bs a^T
								\left(
									\begin{array}{c}
										B_1\\
										B_2
									\end{array}
								\right)(u),	\\					
						\left(
							\begin{array}{c}
								B_3\\
								B_4
							\end{array}
						\right)(u\cdot g)& = 
								\frac{c}{\det(\bs b)}\bs b^T
								\left(
									\begin{array}{c}
										B_3\\
										B_4
									\end{array}
								\right)(u).	
							\end{alignedat}		
					\end{equation*}
	Since $\SL(2,\R)$ acts transitively on $\R^2\backslash\{\bs 0\}$, we can always,
	by choosing an appropriate $0$-adapted coframing,
	reduce to one of the following 4 cases:
		\begin{enumerate}[\qquad\qquad]
			\item[$\mathscr{A}_1$:]{$B_1 = B_3 = 1$, $B_2 = B_4 = 0$;}
			\item[$\mathscr{A}_2$:]{$B_i = 0$ $(i = 1,\ldots,4)$;}
			\item[$\mathscr{A}_3$:]{$B_2, B_3, B_4 = 0$, $B_1 = 1$;}
			\item[$\mathscr{A}_4$:]{$B_1, B_2, B_4 = 0$, $B_3 = 1$.}	
		\end{enumerate}			
	
	Among these cases, $\mathscr{A}_3$ and $\mathscr{A}_4$ are essentially equivalent, because they turn
	into each other as one switches the submersions $\pi$ and $\bar\pi$. Moreover, we will regard $\mathscr{A}_1$
	as the generic case.
	
	We are only interested
	in those rank $2$ B\"acklund transformations that are \emph{not} obtained from a $1$-parameter family of 
	rank $1$ B\"acklund transformations in the obvious way, so we present the following criterion.
	
	\begin{lemma} \label{foliatedrank2}
		Locally a Type $\mathscr{A}$ rank $2$ B\"acklund transformation $(N,\B;\pi,\bar\pi)$
		is one that arises from a $1$-parameter family of rank $1$ 
		B\"acklund transformations (in the sense of Section \ref{familyconstrsection})
		if and only if, for any $0$-adapted coframing \[\bs\omega = (\omega^0,\bar\omega^0,\gamma,\omega^1,\ldots,\omega^4),\] 
			 $\gamma$ is integrable.
	\end{lemma}		 
			
		\emph{Proof}. Suppose that $\gamma\in \Omega^1(U)$ is integrable, 
			that is, there exists a function $\lambda:U\rightarrow \R$ such that $\gamma = f\ed \lambda$
			for some nonvanishing function $f$.
			Let
				\[
				 	\iota_\mu: U_\mu:=\lambda^{-1}(\mu)\hookrightarrow U
				\]
			be the obvious inclusions; define on each $U_\mu$ the differential ideal
				\[
					\B_\mu:= \iota_\mu^*\B.
				\]
			It is straightforward to check that $(U_\mu, \B_\mu; \pi\circ\iota_\mu, \bar\pi\circ\iota_\mu)$
			corresponds to a $1$-parameter family of rank $1$ B\"acklund transformations, from which
			the ideal $\B$ can be recovered by the construction described in Section \ref{familyconstrsection}.

			Conversely, suppose that $(N,\B;\pi,\bar\pi)$ arises from a 
			$1$-parameter family of rank $1$ 
			B\"acklund transformations in the sense of Section \ref{familyconstrsection}. 
			Let $\lambda$ be the parameter. Then we have that $B^1$ is 
			generated by the pullback of $\theta$, $\bar\theta$ and $\ed \lambda$. It follows that, 
			given a $0$-adapted coframing defined on $U\subset N$,
			there exists a linear combination of $\omega^0,\bar\omega^0,\gamma$
			that is integrable. By \eqref{generalomega0}, \eqref{generalbomega0} and \eqref{generaldgamma}, 
			the first derived system of 
			$\lbb\omega^0,\bar\omega^0,\gamma\rbb$ must be $\lbb\gamma\rbb$.
			It is necessary that $\gamma$ is integrable. \qed

\subsection{The case $\mathscr{A}_1$}	
	For the rest of Section \ref{typeAsection}, we will concentrate on the case $\mathscr{A}_1$.
	In this case, $(B_1,B_2) = (B_3, B_4) = (1,0)$.
	In geometric terms, we have reduced to a $G_1$-structure 
	$\mathcal{G}_1\subset \mathcal{G}$, where $G_1$ is the matrix group formed by
	those elements $g\in G$ (see \eqref{trans}) that satisfy 
	 \[
		a_{12} = b_{12} = 0, \qquad a_{11} = b_{11}, \qquad c = a_{22} = b_{22}.
	 \] 
	Thus, $G_1$ is $4$-dimensional. Let $\mathfrak{g}_1$ denote the its Lie algebra. 
	We say that the elements of $\mathcal{G}_1$ are
	\emph{$1$-adapted} coframes of the underlying B\"acklund transformation.
	
	Let $\bs\omega = (\omega^0,\bar\omega^0, \gamma,\omega^1,\ldots,\omega^4)$
	denote the tautological $1$-form on $\mathcal{G}_1$.  The following
	 structure equation holds:
	 \begin{equation}\label{strEq:G1}
	 	\ed \bs\omega = -\bs\varphi \W\bs\omega +  \bs\tau,
	 \end{equation}
	 where $\bs\omega$ should be viewed as a column and 
	\begin{equation}
		\bs\varphi = \left(
					\begin{array}{ccccccc}
						\alpha+\phi &0&0&0&0&0&0\\	
						0&\alpha+\phi&0&0&0&0&0\\	
						0&0&\phi&0&0&0&0\\
						0&0&0&\alpha& 0&0&0\\	
						0&0&0&\beta_1& \phi&0&0\\	
						0&0&0&0&0&\alpha&0\\	
						0&0&0&0&0&\beta_2&\phi	
					\end{array}
				\right),
		\qquad 
		\bs \tau = 
				\left(
					\begin{array}{c}
						\Omega^0\\
						\bar\Omega^0\\
						\Gamma\\
						\Omega^1\\
						\Omega^2\\
						\Omega^3\\
						\Omega^4
					\end{array}
				\right).
		\label{strEq:G1Detail}		
	\end{equation}
			
The $1$-forms $\alpha, \phi, \beta_1,\beta_2\in \Omega^1(\mathcal{G}_1)$
are linearly independent of the components of $\bs\omega$ and among themselves.
The $\mathfrak{g}_1$-valued $1$-form $\bs\varphi$ is called the pseudo-connection $1$-form of 
$\mathcal{G}_1$; $\bs\tau$ is called the torsion.

By the reproducing property of $\bs\omega$
and \eqref{generalomega0}, \eqref{generalbomega0} and \eqref{generaldgamma}, 
the $2$-forms $\Omega^0, \bar\Omega^0$ and $\Gamma$ are restricted. 
By adding any linear combination of
 $\omega^0,\bar\omega^0,\gamma$, $\omega^1,\ldots,\omega^4$ into $\alpha$ and $\phi$,
we can arrange that
\begin{equation}\label{caseA1Tor}
	\left\{\begin{alignedat}{2}
	\Omega^0&=  A\omega^1\W\omega^2 + \omega^3\W(\omega^4 +\gamma) 
		+ (K\gamma+P_i\omega^i)\W\omega^0,\\
	\bar\Omega^0 &= \omega^1\W(\omega^2+\gamma) + \epsilon A \omega^3\W\omega^4 
		-(K\gamma+P_i\omega^i)\W\bar\omega^0,\\
	\Gamma&=C_0\omega^0\W\bar\omega^0 +C_i\omega^i \W\omega^0 
		+D_i \omega^i\W\bar\omega^0,	\qquad  (A>0, A^2\ne \epsilon)
	\end{alignedat}	\right.
\end{equation}
for some functions
$A, C_0, C_i, D_i, K, P_i$ $(i = 1,\ldots,4)$ defined on $\mathcal{G}_1$.  


Further restrictions come from the \emph{invariance property} (see \eqref{InvProp1} and \eqref{InvProp2}),
as the following lemma shows.
 
 \begin{lemma} \label{Rank2GenericInvLemma}
 The following hold
 for the torsion components of $\mathcal{G}_1$.	
 \begin{enumerate}[\rm (A)]
 \item{$\Omega^1$ and $\Omega^2$, when reduced modulo $\{\omega^0,\omega^1,\omega^2\}$, 
 	are both congruent to scalar multiples of $\omega^3\W(\omega^4+\gamma)$;}

\item{$\Omega^1$ 
and $\Omega^2+\Gamma$,
when reduced modulo $\{\bar\omega^0,\omega^1,\omega^2+\gamma\}$,
are both congruent to scalar multiples of $\omega^3\W\omega^4$;}
	
\item{$\Omega^3$ and $\Omega^4$, 
when reduced modulo $\{\bar\omega^0,\omega^3,\omega^4\}$,
are both congruent to scalar multiples of $\omega^1\W(\omega^2+\gamma)$;}

\item{$\Omega^3$ and $\Omega^4 + \Gamma$, 
when reduced modulo $\{\omega^0,\omega^3,\omega^4+\gamma\}$, 
are both congruent to scalar multiples of $\omega^1\W\omega^2$.}
\end{enumerate}
\end{lemma}

\emph{Proof}. We prove (A) and (D), leaving (B) and (C) to the reader, since they are similar.

Let $\sigma:U\rightarrow \mathcal{G}_1$ be any section.
By construction, the vector fields $\bar X_0, X_\gamma$ on $U$ are annihilated by the 
$1$-forms $\sigma^*\omega^0,\sigma^*\omega^1,\sigma^*\omega^2,\sigma^*\omega^3,\sigma^*(\omega^4+\gamma)$ and are dual to 
		$\sigma^*\bar\omega^0$ and $\sigma^*\gamma$.
		Using \eqref{strEq:G1} and dropping the pullback symbols for clarify, we obtain
\begin{align}
	\left\{\begin{alignedat}{1}
		\mathcal{L}_{\bar X_0}\omega^0&\equiv 0\\
		\mathcal{L}_{ X_\gamma}\omega^0&\equiv 0 
	\end{alignedat}\right\}&\mod\omega^0,\label{invLemma1}\\
	\left\{\begin{alignedat}{1}
		\mathcal{L}_{\bar X_0}\omega^1	&\equiv \bar X_0\lefthook\Omega^1\\
		\mathcal{L}_{\bar X_0}\omega^2	&\equiv \bar X_0\lefthook\Omega^2\\
		\mathcal{L}_{X_\gamma}\omega^1	&\equiv  X_\gamma\lefthook\Omega^1\\
		\mathcal{L}_{X_\gamma}\omega^2	&\equiv  X_\gamma\lefthook\Omega^2\\
	\end{alignedat}\right\}&\mod\omega^0,\omega^1,\omega^2,\label{invLemma2}\\
	\left\{\begin{alignedat}{1}				
		\mathcal{L}_{\bar X_0}\omega^3	&\equiv \bar X_0\lefthook\Omega^3\\
		\mathcal{L}_{\bar X_0}(\omega^4+\gamma)&\equiv 
						\bar X_0\lefthook (\Omega^4+\Gamma)\\
		\mathcal{L}_{X_\gamma}\omega^3	&\equiv X_\gamma\lefthook\Omega^3\\
		\mathcal{L}_{X_\gamma}(\omega^4+\gamma)&\equiv 
						X_\gamma \lefthook (\Omega^4+\Gamma)
	\end{alignedat}\right\}&\mod\omega^0,\omega^3,\omega^4+\gamma.\label{invLemma3}
\end{align}

There exist functions $a_i,b_j$ such that 
	\begin{equation*}
		\begin{alignedat}{2}
		\Omega^1& \equiv a_1\bar\omega^0\W\gamma + a_2\bar\omega^0\W\omega^3 + 
							a_3\bar\omega^0\W(\omega^4+\gamma)\\
				& + a_4\gamma\W\omega^3+a_5\gamma\W(\omega^4+\gamma)
						+a_6\omega^3\W(\omega^4+\gamma)	&&\mod\omega^0,\omega^1,\omega^2,
				\\
		\Omega^3& \equiv b_1\bar\omega^0\W\gamma + b_2\bar\omega^0\W\omega^1 + 
							b_3\bar\omega^0\W\omega^2\\
				& + b_4\gamma\W\omega^1+b_5\gamma\W\omega^2
						+b_6\omega^1\W\omega^2
						&&\mod\omega^0,\omega^3,\omega^4+\gamma.				
		\end{alignedat}					
	\end{equation*}

From this we obtain
	\begin{align*}
	&\left.
	\begin{alignedat}{2}
		\bar X_0\lefthook\Omega^1&\equiv a_1\gamma+a_2\omega^3+a_3(\omega^4+\gamma)\\
		X_\gamma\lefthook\Omega^1&\equiv - a_1\bar\omega^0 
							+a_4\omega^3+a_5(\omega^4+\gamma)
	\end{alignedat}
	\right\}&&\mod\omega^0, \omega^1,\omega^2.	
	\\
	&\left.
	\begin{alignedat}{2}
		\bar X_0\lefthook\Omega^3&\equiv b_1\gamma+b_2\omega^1+b_3\omega^2\\
		X_\gamma\lefthook\Omega^3&\equiv - b_1\bar\omega^0 
							+b_4\omega^1+b_5\omega^2
	\end{alignedat}
		\right\}&&\mod\omega^0, \omega^3,\omega^4+\gamma.	
	\end{align*}	
	
In order for the right-hand-sides of \eqref{invLemma2} and \eqref{invLemma3} to be congruent to zero modulo 
$\{\omega^0,\omega^1,\omega^2\}$ and $\{\omega^0,\omega^3,\omega^4+\gamma\}$, 
respectively,  it is necessary that \[a_1 = \cdots = a_5 = b_1 = \cdots = b_5 = 0.\]
This justifies (A) for $\Omega^1$ and (D) for $\Omega^3$; the argument for $\Omega^2$ and $\Omega^4+\Gamma$
is similar.
	\qed\\

We remark that the converse of Lemma \ref{Rank2GenericInvLemma} is also true, in the sense that
any $G_1$-structure on a $7$-manifold whose tautological $1$-form $\bs\omega$ satisfies
\eqref{strEq:G1}, \eqref{caseA1Tor} and the conditions (A)-(D) in Lemma \ref{Rank2GenericInvLemma}
is associated to a Type $\mathscr{A}_1$ B\"acklund transformation.

By Lemma \ref{Rank2GenericInvLemma}, on $\mathcal{G}_1$ there exist functions 
\[	
	\begin{alignedat}{2}
		&T^i_0, S_i\qquad \qquad 
			&&(i = 1,\ldots,4) \\
		T^\sigma_\rho, &\bar T^\sigma_\rho, R^\sigma_\rho 
			&&(\sigma,\rho\in \{1,2\} \text{ or } \{3,4\})\\
		&T^i_{jk}
			&&(i,j,k = 1,\ldots,4;~ T^i_{jk} + T^i_{kj} = 0)
	\end{alignedat}
\]
such that $\Omega^1,\ldots,\Omega^4$ take the form:
\begin{equation}\label{caseA1Om1to4}{\small
	\left\{\;
	\begin{alignedat}{2}
		\Omega^1&= T^1_0\omega^0\W\bar\omega^0 + \omega^0\W(T_1^1 \omega^1+T_2^1(\omega^2+\gamma)) + \bar\omega^0\W(\bar T_1^1\omega^1+\bar T_2^1\omega^2) \\	
		& +\gamma\W(R_1^1\omega^1+R_2^1\omega^2)
		+ S_1\omega^3\W(\gamma+\omega^2+\omega^4)
		+ \frac{1}{2}T^1_{ij}\omega^i\W\omega^j,\\
		\Omega^2&= T^2_0\omega^0\W\bar\omega^0 +\omega^0\W(C_3\omega^3+C_4\omega^4+T_1^2\omega^1+T_2^2(\omega^2+\gamma))\\
		&+\bar\omega^0\W(\bar T_1^2\omega^1+\bar T_2^2\omega^2)
				 +\gamma\W(C_2\omega^0 + R_1^2\omega^1+R_2^2\omega^2) \\
			&+S_2\omega^3\W(\gamma+\omega^2+\omega^4)	 
			+\frac{1}{2}T_{ij}^2\omega^i\W\omega^j,\\			
		\Omega^3&= T^3_0 \bar\omega^0\W\omega^0 +\bar\omega^0\W(\bar T_3^3\omega^3+\bar T_4^3(\omega^4+\gamma))+ \omega^0\W(T_3^3\omega^3+T_4^3\omega^4)\\
	&+\gamma\W(R_3^3\omega^3+R_4^3\omega^4)
	+S_3\omega^1\W(\gamma+\omega^2+\omega^4)
	+\frac{1}{2}T_{ij}^3\omega^i\W\omega^j,\\
\Omega^4&=	T^4_0 \bar\omega^0\W\omega^0 
	+\bar\omega^0\W(D_1\omega^1+D_2\omega^2
	+\bar T_3^4\omega^3 +\bar T_4^4(\omega^4+\gamma))\\
		&+ \omega^0\W(T_3^4\omega^3+T_4^4\omega^4) 
		 +\gamma\W(D_4\bar\omega^0+R_3^4\omega^3+R_4^4\omega^4)\\
		 &+S_4\omega^1\W(\gamma+\omega^2+\omega^4)
		  +\frac{1}{2}T_{ij}^4\omega^i\W\omega^j,
\end{alignedat}\right.}
\end{equation}		

where
\[
T^1_{23}, T^1_{24}, T^1_{34}, T^2_{23}, T^2_{24}, T^2_{34}, T^3_{12}, T^3_{14}, T^3_{24}, T^4_{12}, T^4_{14}, T^4_{24}\] 
are all zero. 

The coefficients in $\Omega^1,\ldots,\Omega^4$ are not all determined.
In fact, by adding a linear combination of 
$\omega^0,\bar\omega^0,\gamma,\omega^1,\ldots,\omega^4$ into $\beta_1$, 
we can arrange that
		\[
			T^2_{12} = T^2_{13} =  T^2_{14} =  T^2_1 =  \bar T^2_1 = R^2_1 = 0;
		\] 
 in a similar manner, by adjusting $\beta_2$, we can arrange that  
		\[
			  T^4_{13} =  T^4_{23} =  T^4_{34} =\bar T^4_3 = T^4_3 =  R^4_3 = 0;
		\]
by adding a multiple of $\gamma$ into $\alpha$ and subtracting the same multiple from 
$\phi$, we can arrange that
		\[
			R^4_4 = - R^2_2.
		\]

As a standard step in the method of equivalence, 
we apply the \emph{fundamental identity} $\ed(\ed \bs\omega) = 0$ to \eqref{strEq:G1}. From this, we obtain relations 
among the \emph{torsion functions}, that is, 
the coefficient functions occurring in $\Omega^0, \bar\Omega^0,\Gamma,\ldots,\Omega^4$.

In particular, by expanding the expressions
\begin{equation*}\left\{
	\begin{alignedat}{2}
	\ed(\ed\omega^0)&\mod \omega^0,\omega^1,\\ 
		\ed(\ed\omega^0)&\mod\omega^0,\omega^2,\\
	\ed(\ed\bar\omega^0)&\mod\bar\omega^0,\omega^3,\\ 
		\ed(\ed\bar\omega^0)&\mod\bar\omega^0,\omega^4,\\
	\ed(\ed\gamma) &\mod \omega^0,\bar\omega^0,\gamma,
	\end{alignedat}\right.
\end{equation*}
we find the relations
\begin{equation}\label{caseA1CandDrelate}
			D_2 = -\frac{\epsilon}{A} C_2, \quad
			 C_4 = -\frac{1}{A} D_4, \quad
			  D_1 = -\frac{\epsilon}{A} C_1, \quad
			   C_3 = -\frac{1}{A} D_3
\end{equation}
and
	\begin{equation}\label{caseA1Tsolved}
		\left\{\;
		\begin{alignedat}{2}	
			T^1_{12}& = P_2-K-R^1_1 - R^2_2, \\
			T^3_{34} &= -P_4+K-R^3_3+R^2_2,\\
			T^1_{13}& = -\epsilon AS_4+ S_2+P_3,\\
			T^3_{13}& = AS_2 - S_4+P_1,\\
			T^1_{14}& = \epsilon A S_3+P_4,\\
			T^3_{23}& = -AS_1+P_2,\\
			 T^2_2 &= C_2 - T^1_1, \\
			 \bar T^4_4& = D_4-\bar T^3_3 , \\
		\end{alignedat}
		\right.
	\end{equation}

Thus, we replace the functions on the left-hand-side of \eqref{caseA1CandDrelate} and \eqref{caseA1Tsolved} by 
the corresponding expressions on the right-hand-side.
	
Next, by expanding the expressions
		\[
			\left\{
			\begin{alignedat}{1}
				\ed(\ed\omega^0) &\mod\omega^0, \\
				 \ed(\ed\bar\omega^0) &\mod \bar\omega^0, 
			\end{alignedat}
			\right.
		\]
 we obtain
\begin{equation}
	R^2_2 = K+\frac{1}{2}(R^3_3 - R^1_1) + \frac{1}{2A}(\epsilon S_1 - S_3)
\end{equation}
 and
	\begin{equation}\label{caseA1dA}\left\{
		\begin{alignedat}{2}
			\ed A&= -A(T^3_3+T^4_4)\omega^0 -A(\bar T^1_1+\bar T^2_2)\bar\omega^0
			\\
			& -\frac{1}{2}(AR^1_1+AR^3_3 +\epsilon S_1+S_3)\gamma
			\\
			& - (2AP_1 +(A^2- \epsilon)S_2)\omega^1 
			 - (2AP_2 - (A^2 - \epsilon)S_1)\omega^2
			 \\
			 & +(2AP_3- (\epsilon A^2 - 1)S_4)\omega^3
			 +(2AP_4+(\epsilon A^2 - 1)S_3)\omega^4.
		\end{alignedat}\right.
	\end{equation}	
Therefore, $A$ is a scalar invariant of the $G_1$-structure on $N$. 

\begin{lemma}\label{CaseAC0lemma}
	In \eqref{caseA1Tor}, 
	\begin{equation}\label{C0inothers}
		C_0 = \frac{\epsilon C_1S_1+\epsilon C_2S_2 - D_3S_3 - D_4S_4}{A}.
	\end{equation}
\end{lemma}
\emph{Proof}.
Computing $\ed(\ed\gamma)$ and reducing in two ways, 
modulo $\{\bar\omega^0,\gamma,\omega^1,\omega^2\}$ and 
$\{\omega^0,\gamma,\omega^1,\omega^2\}$, respectively,
we obtain
	\begin{equation*}\left.
		\begin{alignedat}{2}
		\omega^3\W\ed D_3 + \omega^4\W\ed D_4
			&\equiv (f_1+g)\omega^3\W\omega^4 
			+ \psi
			\\
		\omega^3\W\ed D_3 + \omega^4\W\ed D_4
			&\equiv(f_2+g)\omega^3\W\omega^4 
			+ \psi 
		\end{alignedat}\right\} \mod \omega^0,\bar\omega^0,\omega^1,\omega^2,\gamma,
	\end{equation*}
where
	\begin{equation*}\left\{
		\begin{alignedat}{2}
		g& = 2D_3K - \frac{D_3}{2}(R^1_1+R^3_3) - D_4P_3 - \frac{D_3}{2A}(S_3 - \epsilon S_1),\\
		f_1&=\epsilon A^2C_0-A(C_1S_1 + C_2S_2) + \frac{\epsilon A^2 - 1}{A}(D_3S_3+D_4S_4),\\
		f_2&= C_0 - \frac{\epsilon}{A}(C_1S_1+C_2S_2),
		\end{alignedat}\right.
	\end{equation*}
and 
	\begin{equation*}
		\psi= -(2D_3\alpha+D_4\beta_2)\W\omega^3	 - D_4(\alpha+\phi)\W\omega^4.	
	\end{equation*}	
Evidently, $f_1$ and $f_2$ must be equal. This implies \eqref{C0inothers}, because $A^2\ne \epsilon$.		
\qed	
\\

Now, other than $A$, there are $31$ torsion functions remaining.  By applying the fundamental identity to 
\eqref{strEq:G1}, we find the infinitesimal transformation of those functions
under the action of $G_1$, listed as follows, where all congruences are reduced modulo the $1$-forms $\omega^0,
\bar\omega^0,\gamma,\omega^1,\ldots,\omega^4$.
\begin{equation}\label{A1torsiontrans}
		\left\{\;
		\begin{alignedat}{5}
		\ed K&\equiv \phi K,\\
		\ed S_1&\equiv \phi S_1,\qquad &&\ed S_2&&\equiv -\beta_1S_1 +\alpha S_2,\\
		\ed S_3&\equiv \phi S_3,\qquad &&\ed S_4&&\equiv -\beta_2 S_3+\alpha S_4,\\
		\ed C_1&\equiv 2 \alpha C_1+\beta_1 C_2,\qquad &&\ed C_2&&\equiv (\alpha+\phi) C_2,\\
		\ed D_3&\equiv 2\alpha D_3 +\beta_2 D_4, \qquad &&\ed D_4&&\equiv (\alpha+\phi) D_4,\\
		\ed P_1&\equiv \alpha P_1 + \beta_1 P_2,\qquad&&\ed P_2&&\equiv \phi P_2,\\
		\ed P_3&\equiv \alpha P_3+\beta_2P_4,\qquad &&\ed P_4&&\equiv \phi P_4,\\
		\ed R^1_1&\equiv \phi R^1_1+\frac{\beta_1}{2}R^1_2 - \frac{\beta_2}{2}R^3_4,\qquad 
					&&\ed R^1_2&&\equiv (2\phi - \alpha) R^1_2,\\
		\ed R^3_3&\equiv \phi R^3_3 - \frac{\beta_1}{2}R^1_2 +\frac{\beta_2}{2}R^3_4,\qquad
					&&\ed R^3_4&&\equiv(2\phi - \alpha) R^3_4,\\
		\ed T^1_0&\equiv (\alpha +2\phi) T^1_0, \qquad &&\ed T^2_0&&\equiv 
					(2\alpha+\phi)T^2_0-\beta_1 T^1_0,\\
		\ed T^3_0&\equiv (\alpha+2\phi) T^3_0, \qquad &&\ed T^4_0&&\equiv 
					(2\alpha+\phi)T^4_0-\beta_2 T^3_0,	\\
		\ed T^1_1&\equiv (\alpha+\phi) T^1_1+\beta_1 T^1_2,\qquad
					&&\ed T^1_2&&\equiv  2\phi T^1_2,\\
		\ed \bar T^1_1&\equiv (\alpha+\phi) \bar T^1_1+\beta_1\bar T^1_2,\qquad
					&&\ed \bar T^1_2 &&\equiv 2\phi \bar T^1_2,\\
		\ed \bar T^2_2&\equiv (\alpha+\phi)\bar T^2_2 - \beta_1\bar T^1_2,\qquad
					&&\ed T^4_4&&\equiv (\alpha+\phi)T^4_4 - \beta_2T^3_4,\\				
					\ed T^3_3&\equiv (\alpha+\phi) T^3_3+\beta_2 T^3_4,\qquad
					&&\ed T^3_4&&\equiv 2\phi T^3_4,\\			
		\ed \bar T^3_3&\equiv (\alpha+\phi) \bar T^3_3+\beta_2\bar T^3_4,\qquad
				&&\ed \bar T^3_4&&\equiv	2\phi \bar T^3_4.					
	\end{alignedat}\right.
\end{equation}

These $31$ functions are called the \emph{relative invariants}
of $\mathcal{G}_1$.

It is natural to impose $G_1$-invariant algebraic constraints on these 
relative invariants and then study the resulting structures.

\subsection{Existence and rigidity}
By \eqref{A1torsiontrans}, $C_2, D_4$ and $K$ scale under the $G_1$-action along each fibre of
$\mathcal{G}_1$.
It follows that the open condition:
\begin{equation}\label{C2D4K}
C_2, D_4, K\ne 0
\end{equation}
is $G_1$-invariant. Assuming this and using \eqref{A1torsiontrans}, we can normalize to 
\[
D_4 = K = 1, \qquad C_1 = C_3 = 0.
\]
This normalization exhausts the $G_1$-action on each fibre of $\mathcal{G}_1$;
thus, it gives rise to a canonical coframing $\bs\omega$ defined on $N$. 
The original pseudo-connection
components $\alpha,\phi,\beta_1,\beta_2$ are now linear combinations of the components of $\bs\omega$.
All coefficients that occur in $\ed\bs\omega$ are then explicit expressions of $54$
\emph{primary invariants} (including $A, C_2, S_i$, etc.) of the underlying B\"acklund transformation. 
Once we take into account the fundamental identity $\ed(\ed\bs\omega) = 0$, 
relations among these primary invariants arise. At a stage, we find $7$ such relations. For the remaining $47$
primary invariants, say, $X_i$ $(i = 1,\ldots,47)$, we find a condition on their covariant derivatives 
that is both necessary and sufficient for $\ed(\ed\bs\omega)$ to vanish identically. 
However, the computation of $\ed(\ed X_i)$ becomes difficult to manage even with the help of a machine. 
Thus, for the time being  it seems unlikely to achieve a generality result under the only assumption \eqref{C2D4K} by directly 
applying Cartan's generalization of Lie's third theorem, but one can still investigate 
under less generic assumptions.

Still assuming \eqref{C2D4K} and in view of Lemma \ref{CaseAC0lemma}, 
a simple constraint that we can impose is 
	\begin{equation}\label{SiAssumption}
		S_i= 0 \quad (i = 1,\ldots, 4).
	\end{equation}
By Lemma \ref{foliatedrank2}, the B\"acklund transformations, if any, that belong to this case are \emph{not}
obtained in the obvious way from a $1$-parameter family of rank $1$ B\"acklund transformations.
In fact, those B\"acklund 
transformations that satisfy \eqref{C2D4K} and \eqref{SiAssumption}
not only exist but also display some rigidity, as the following theorem indicates.

\begin{theorem}\label{caseA1rigidityThm}
	Up to contact equivalence, 
	 the space of local
	Type $\mathscr{A}_1$ B\"acklund transformations whose relative invariants  satisfy
	\begin{equation}\label{caseA1ThmAssu}
		\left\{
			\begin{alignedat}{2}
				&C_2, D_4, K\ne 0,\\
				&	S_i=0 \quad (i = 1,\ldots, 4),
			\end{alignedat}
			\right.		
	\end{equation}	
	is parametrized by a finite number of constants.
\end{theorem}

\emph{Proof}. By the assumption, 
we can apply a pointwise $G_1$-action to reduce to the subbundle 
$\mathcal{G}_2\subset \mathcal{G}_1$ defined by
\begin{equation*}
     \left\{\begin{alignedat}{2}
     			D_4 &= K = 1,\\ 
			 C_1 &= D_3 = 0.
		\end{alignedat}
	\right.		 
\end{equation*}
The result is an $\{e\}$-structure. On $\mathcal{G}_2$, which is diffeomorphic to $N$,
the $1$-forms $\phi, \alpha,\beta_1, \beta_2$ can be expressed as
\begin{equation}\label{caseA1Thmphab1b2}
	\left\{\begin{alignedat}{2}
		\phi& = M_0\omega^0+\bar M_0\bar\omega^0
				+ M_\gamma\gamma+ M_i\omega^i,\\
		\alpha& = U_0\omega^0+\bar U_0\bar\omega^0
				+ U_\gamma\gamma+ U_i\omega^i,\\	
		\beta_1& = V_0\omega^0+\bar V_0\bar\omega^0
				+ V_\gamma\gamma+ V_i\omega^i,\\	
		\beta_2& = W_0\omega^0+\bar W_0\bar\omega^0
				+ W_\gamma\gamma+ W_i\omega^i,\\				
	\end{alignedat}\right.
\end{equation}
for certain coefficient functions defined on $N$.
We note that $V_1, W_3$ do not appear in the structure equations. There are 
$50$ primary invariants for the $\{e\}$-structure:
$26$ from the coefficients of \eqref{caseA1Thmphab1b2};
$23$ from what's remaining in \eqref{A1torsiontrans};
and the function $A$.
For convenience, we rename them as $X_1,\ldots, X_{50}$ with the particular assignment:
	\begin{equation}\label{X1to4}
		X_1:= A, \qquad X_2 := C_2, \qquad X_3 := M_0, \qquad X_4 := \bar M_0.
	\end{equation}
Furthermore, define the covariant derivatives $X_{i,j}$ by
	\begin{equation*}
		\begin{alignedat}{2}
		\ed X_\rho = X_{\rho,1}\omega^0+ X_{\rho,2}\bar \omega^0
				+ X_{\rho,3}\gamma + X_{\rho,4}\omega^1
				+X_{\rho,5}\omega^2+&X_{\rho,6}\omega^3
				+X_{\rho,7}\omega^4 \\
				& (\rho = 1,\ldots,50).
		\end{alignedat}		
	\end{equation*}

Expanding the equation 
\[
	\ed(\ed\bs\omega) = 0
\]
yields a system of polynomial equations in the $X_\rho$ and $X_{\rho,i}$.	
This system imposes relations among the primary invariants
$X_\rho$; we update these relations and repeat until a stage when 
the system of polynomial equations arising at that stage can be 
solved for the covariant derivatives $X_{\rho, i}$
that appear in the system.

Following this procedure, after a somewhat lengthy calculation
using Maple\texttrademark, we find that (see also \eqref{caseA1oms} and 
\eqref{caseA1Hs} below)
\begin{itemize}
\item{all $X_\rho$ $(\rho = 1,\ldots,50)$
are explicit expressions of $X_1,\ldots,X_4$;}
\item{all covariant derivatives of 
		$X_1,\ldots, X_4$ are explicit expressions of $X_1,\ldots, X_4$;}
\item{$\ed(\ed\bs\omega)$ and $\ed(\ed X_i)$ $(i = 1,\ldots,4)$ vanish identically.}		
\end{itemize}

Since the systems of structure equations is involutive, and 
since there are no free derivatives of $X_1,\ldots, X_4$, the structure of 
the underlying B\"acklund transformation near $p\in N$ is completely determined by 
the value of $(X_i)$ $(i = 1,\ldots, 4)$ at $p$. This is a consequence of 
Cartan's third fundamental theorem (see \cite{br14}). The conclusion follows.
\qed
\\

On the next two pages we record the structure equations \eqref{caseA1oms} and 
\eqref{caseA1Hs} that are mentioned 
in the proof of Theorem \ref{caseA1rigidityThm}. The substitutions
\begin{equation}\label{subHtoX}
	H_1:= X_1, \quad H_2:= X_1X_2, \quad H_3:= X_1X_3, \quad H_4:= \frac{X_1X_4}{X_2}
\end{equation}
are made to simplify the expressions.

\newpage
\vspace*{\fill}
\begin{equation}\label{caseA1oms}
{\footnotesize\left\{\;
\begin{alignedat}{2}
	\ed \omega^0
		&= \frac{1}{2}\(1 - \frac{\epsilon H_2}{(H_1)^2}\)	
			\omega^0\W\bar\omega^0 + (3\epsilon H_4 - H_3+1)\omega^0\W\gamma	
			+ H_1\omega^1\W\omega^2
		\\
		&- \frac{1}{2}\(H_2(\epsilon H_4 - 3H_3 - 3) - 3\epsilon H_4+ H_3 - 3\)\omega^0\W\omega^2
		+ \omega^3\W(\omega^4+\gamma)\\
		&-\frac{1}{2}\(H_3 -3\epsilon H_4  +1
				+\frac{(H_1)^2\(\epsilon H_3- 3 H_4 - 3\epsilon \)}{H_2}\)\omega^0\W\omega^4,\\[0.8em]
\ed\bar\omega^0&= 	(3\epsilon H_4 - H_3+3)\bar\omega^0\W(\omega^2+\gamma)
					+ (\epsilon H_4 + H_3 +1)\bar\omega^0\W\omega^4\\
				& +\omega^1\W(\omega^2+\gamma)
				+\epsilon H_1\omega^3\W\omega^4,\\[0.8em]
\ed\gamma& =  \frac{1}{H_1}(H_3\gamma +H_2\omega^2 - \omega^4)\W\omega^0						
			 +\(\frac{H_2H_4}{(H_1)^2}\gamma
			 		 - \frac{\epsilon H_2}{(H_1)^2}\omega^2 +\omega^4\)
					 \W\bar\omega^0
		\\
		 &+\(\frac{3}{4} +\frac{\epsilon H_2}{4(H_1)^2}\)\omega^1\W\gamma		  
		 	+(3\epsilon H_4 - H_3+3)(\epsilon H_4+1) \omega^2\W\gamma
		 \\
		 & -\frac{3H_2+1}{4H_1}\omega^3\W\gamma
		  +(\epsilon H_4 - 3H_3 - 3)(H_3+1)\omega^4\W\gamma,\\[0.8em]
\ed\omega^1&=
			 (3(H_4)^2 - 3(H_3)^2 + 8\epsilon H_4 - 6H_3+1)\omega^1\W\gamma	 
			+ \frac{H_3+1}{H_1}\omega^0\W\omega^1	 \\
			&+(3(H_4)^2 - \epsilon H_3H_4 + 7\epsilon H_4 +4)\omega^1\W\omega^2
			+\frac{ H_2( H_4+\epsilon)}{(H_1)^2}\bar\omega^0\W\omega^1\\
			&-( 3(H_3)^2  - \epsilon H_3H_4 + 5H_3 - 2\epsilon H_4+2)\omega^1\W\omega^4
			-\frac{3H_2+1}{4H_1}\omega^1\W\omega^3,\\[0.8em]
\ed\omega^2&= -\frac{1}{H_1}\omega^0\W(\omega^4+\gamma)
				+\frac{H_2 - H_3 - 1}{H_1}\omega^0\W\omega^2		
				-\frac{H_2(H_4+\epsilon)}{(H_1)^2}\bar\omega^0\W\omega^2\\
			&+ \(\frac{3}{4} +\frac{\epsilon H_2}{4(H_1)^2}\)\omega^1\W\omega^2
			+ (3(H_3)^2 - 3(H_4)^2 + 7 H_3 - 7\epsilon H_4)\omega^2\W\gamma\\
			&+\frac{3H_2+1}{4H_1}\omega^2\W\omega^3
			 - (\epsilon H_4 - 3H_3 - 3)(H_3+1)\omega^2\W\omega^4,\\[0.8em]
\ed\omega^3&= \frac{1+H_2+2H_3}{2H_1}\omega^0\W\omega^3
				+\(-\frac{1}{2}+\frac{H_2H_4}{(H_1)^2} +\frac{3\epsilon H_2}{2(H_1)^2}\)
				\bar\omega^0\W\omega^3
				\\
			& - \(\frac{3}{4}+\frac{\epsilon H_2}{4(H_1)^2}\)\omega^1\W\omega^3	
			+(3(H_4)^2  + 8\epsilon H_4- 3(H_3)^2 - 6H_3 +3)\omega^3\W\gamma		 
			\\
			&+\frac{1}{2}\(\epsilon H_4 (H_2+ 2H_3 -6\epsilon H_4- 15)
					- 3 H_2H_3  - 3H_2+3H_3 - 9\)\omega^2\W\omega^3\\
			& - \frac{1}{2}\(\frac{(H_1)^2(\epsilon H_3 -3H_4- 3\epsilon)}{H_2}
		+(3H_3 - \epsilon H_4+3)(2H_3+1)\)\omega^3\W\omega^4,\\[0.8em]
\ed\omega^4&= -\frac{H_3+1}{H_1}\omega^0\W\omega^4
				 - \frac{\epsilon H_2}{(H_1)^2}\bar\omega^0\W(\omega^2+\gamma) 
				 + \(1  - \frac{H_2(H_4+ \epsilon)}{(H_1)^2}\)\bar\omega^0\W\omega^4	\\
		&+\(\frac{3}{4} +\frac{\epsilon H_2}{4(H_1)^2}\)\omega^1\W\omega^4
		+(H_4+\epsilon)(3H_4 - \epsilon H_3 + 3\epsilon)\omega^2\W\omega^4\\
		&	 -\frac{3H_2+1}{4H_1}\omega^3\W\omega^4
		+(3(H_3)^2 - 3(H_4)^2 + 7H_3 - 7\epsilon H_4)\omega^4\W\gamma.
\end{alignedat}\right.}
\end{equation}
\vspace*{\fill}

\newpage
\begin{equation}\label{caseA1Hs}
{\footnotesize\left\{\;
\begin{alignedat}{2}
\ed H_1&= -\frac{1}{2}(H_2-1)\omega^0 
			- \frac{1}{2}\(H_1- \frac{\epsilon H_2}{H_1}\)\bar\omega^0
			+  2H_1(H_3 - \epsilon H_4)\gamma
			\\		
		&  - \frac{H_1}{2}(\epsilon H_2H_4+3\epsilon H_4 - 3H_2H_3 - 3H_2 - H_3+3)
		  	\omega^2\\
		&-\frac{H_1}{2}\(\frac{(H_1)^2}{H_2}(3H_4 +3\epsilon- \epsilon H_3)
			+(\epsilon H_4 - 3H_3 - 3)\)\omega^4,	\\[0.8em]
\ed  H_2&= \frac{H_2(1-H_2)}{H_1}\omega^0
				 + 4H_2(H_3+1)(\omega^4+\gamma)\\
		&	-H_2(\epsilon H_2H_4 - \epsilon H_4 - 3H_2H_3 - 3H_2 - H_3 - 1)\omega^2,	\\[0.8em]
\ed H_3& = \frac{H_3(H_3+1)}{H_1}\omega^0+\frac{H_2(H_3+1)(H_4+\epsilon)}{(H_1)^2}\bar\omega^0	\\
	& - (6\epsilon H_4- 3(H_3)^2 + 3(H_4)^2 - 4H_3+3)(H_3+1)\gamma\\
	&-(H_3+1)\(\frac{3}{4}+ \frac{\epsilon H_2}{4(H_1)^2}\)\omega^1
	+ \frac{3H_2H_3+4H_2+H_3}{4H_1}\omega^3\\
	&-(H_3+1)(3H_4 - \epsilon H_3+3\epsilon)(H_4+\epsilon)\omega^2\\
	& - H_3(H_3+1)(\epsilon H_4 - 3H_3 - 3)\omega^4,\\[0.8em]
\ed H_4& = \frac{(H_3+1)(H_4+\epsilon)}{H_1}\omega^0
			+\frac{H_2H_4(H_4+\epsilon)}{(H_1)^2}\bar\omega^0\\
	& -(4\epsilon H_4 - 3(H_3)^2 + 3(H_4)^2 - 6 H_3 - 3)(H_4+\epsilon)\gamma\\
	&-\(\frac{3H_4}{4} +\frac{\epsilon H_2H_4}{4(H_1)^2}+\epsilon\)\omega^1	
	 - H_4(3\epsilon H_4 - H_3+3)(\epsilon H_4+1)\omega^2\\
	 &+\frac{(H_4+\epsilon)(3H_2+1)}{4H_1}\omega^3
	  - (H_3+1)(H_4+\epsilon)(\epsilon H_4 - 3H_3 - 3)\omega^4.
\end{alignedat}\right.}
\end{equation}

\subsection{Symmetries}\label{TypeA1symmetries}

For each choice of $\epsilon$, 
	\[
		\Theta:=(H_1, H_2, H_3, H_4)
	\]
is a map defined on $N$, that is,
	\[
		\Theta: N\rightarrow \mathcal{S}
						:= \{(h_i)\in \R^4: h_1>0, (h_1)^2\ne \epsilon, h_2\ne 0\}\subset \R^4.
	\]

Generically, $\Theta$ has 
rank $4$. On the other hand, there exist singular loci in 
$\mathcal{S}$ where $\Theta$ drops rank; this corresponds to B\"acklund
transformations with lower cohomogeneity. Now we will describe these cases
in more detail.

Let $\bs h = (h_{i\lambda})$ $(i = 1,\ldots,4;\lambda = 1,\ldots, 7)$ be the coefficient matrix 
occurring in \eqref{caseA1Hs}, whose entries are defined by
	\[
		\ed H_i = h_{i1}\omega^0+ h_{i2}\bar\omega^0
				+h_{i3}\gamma + \sum_{j = 1}^4h_{i(3+j)}\omega^j \qquad(i = 1,\ldots,4).
	\]
	
For $i,j,k,\ell\in\{1,2,\ldots, 7\}$ all distinct, let 
$\hat{\bs h}_{ijk\ell}$ denote the square matrix formed by columns $i,j,k$ and $\ell$ of $\bs h$.
We obtain
	\begin{equation}\label{keyMinor}
		\det(\hat{\bs h}_{1246}) = \frac{H_2\cdot\chi_1\cdot\chi_2
						\cdot\chi_3}{{32 (H_1)^5} },
	\end{equation}
where
	\begin{equation}\label{expressionChi}
		\left\{
		\begin{alignedat}{2}
		\chi_1& = H_2 - 1,\\
		\chi_2&= (H_1)^2 - \epsilon H_2,\\
		\chi_3& = -\epsilon(3H_2H_3 +H_3+4 H_2)\cdot\chi_2\\
						& - (\epsilon H_2 + 3(H_1)^2)(H_4+\epsilon)\cdot\chi_1.
		\end{alignedat}
		\right.			  
	\end{equation}
Since $H_2\ne 0$, in order for $\det(\hat{\bs h}_{1246})$ to vanish, one of $\chi_1, \chi_2, \chi_3$
must vanish. Also note that $\chi_1, \chi_2$ cannot vanish simultaneously.

If $\chi_1 = 0$, that is, $H_2 = 1$, 
then
	\[
		\det( \hat{\bs h}_{2346} )=  - \frac{\epsilon(H_3+1)^2((H_1)^2-\epsilon)^2}{2(H_1)^4}.
	\]
Since $(H_1)^2 \ne \epsilon$, in order for $\rank(\bs h)$ to be less than $4$, it is necessary
that 
	\[
		H_3 = -1.
	\]
One can verify that, wherever $H_2 = 1$ and $H_3 = -1$, the exterior derivatives $\ed H_2$ and 
$\ed H_3$ vanish identically. Moreover, in this case, the matrix $\bs h$
has constant rank $2$ (in particular, the rank is never strictly less than $2$), 
as a result of the assumption $(H_1)^2\ne \epsilon$.
Therefore, the corresponding B\"acklund transformations have cohomogeneity $2$.

If $\chi_2 = 0$, that is, $H_2 = \epsilon (H_1)^2$, then
	\[
		\det(\hat{\bs h}_{1346}) = -\frac{1}{2} ((H_1)^2 - \epsilon)^2H_1(H_4+\epsilon)^2.
	\]
Again, in order for $\rank(\bs h)<4$, it is necessary that 
	\[
		H_4 = -\epsilon.
	\]
By \eqref{caseA1Hs}, wherever $\chi_2 = 0$ and $H_4 = -\epsilon$, the exterior derivatives 
$\ed \chi_2$ and $\ed H_4$ vanish
identically. Moreover, $\rank(\bs h) = 2$, since $(H_1)^2\ne \epsilon$. The corresponding B\"acklund
transformations have cohomogeneity $2$.

If $\chi_3 = 0$ and $\chi_1,\chi_2\ne 0$, then it is straightforward to verify that 
	$\ed \chi_3$ vanishes automatically. We note that $\chi_3 = 0$ includes the cases
	\begin{equation*}
		\left\{
			\begin{alignedat}{2}
				\chi_1 &= 0\\
				H_3 &= -1
			\end{alignedat}
		\right.
		\qquad \text{ and }\qquad 
		\left\{
			\begin{alignedat}{2}
				\chi_2 &= 0\\
				H_4 &= -\epsilon
			\end{alignedat}
		\right.
	\end{equation*}
	above. Hence, the remaining question is whether the image of 
	$\Theta$ can belong to
	any other $k$-dimensional locus $(k<3)$ contained in the zero set of $\chi_3$. 
	To answer this question, let $\hat{\bs h}^{abc}_{ijk}$ denote the matrix
	formed by rows $a,b,c$ and columns $i,j,k$ of $\bs h$. We have
	\begin{equation}\left\{
		\label{caseA1rank3dets}
		\begin{alignedat}{2}	
			\det(\hat{\bs h}^{124}_{236}) &= -\frac{\chi_2H_2\cdot (3H_2+1)(H_3+1)(H_4+\epsilon)}{2(H_1)^2},\\
			\det(\hat{\bs h}^{123}_{126}) &= -\frac{\chi_1\chi_2H_2\cdot((3H_2+1)(H_3+1)+ \chi_1)}{8(H_1)^3},\\
			\det(\hat{\bs h}^{124}_{234}) &= \frac{\epsilon\chi_2 H_2(H_3+1)
							\cdot\((3\epsilon (H_1)^2 + H_2)(H_4+\epsilon)+\chi_2\)}{2(H_1)^3}.
		\end{alignedat}\right.	
	\end{equation}
	Since $\chi_1, \chi_2, H_2$ are nonzero by assumption, it is impossible for the $3$ 
	determinants in \eqref{caseA1rank3dets}
	to vanish simultaneously. To sum up, we have proved the following result.
	
\begin{theorem}\label{caseA1cohomogCls}	
	Let $(N,\B;\pi,\bar\pi)$ be a Type $\mathscr{A}_1$ rank $2$ B\"acklund transformation whose relative 
	invariants satisfy \eqref{caseA1ThmAssu}, the canonical structure equations on $N$
	 being \eqref{caseA1oms} and \eqref{caseA1Hs}. This B\"acklund transformation is of
	 \begin{itemize}
	 	\item{cohomogeneity $2$ if and only if
		\begin{equation}\label{caseA1cohomog2Cond}
			\text{either}\quad
		\left\{
			\begin{alignedat}{2}
				H_2 &= 1\\
				H_3 &= -1
			\end{alignedat}
		\right.
		\quad \text{ or }\quad 
		\left\{
			\begin{alignedat}{2}
				H_2 &= \epsilon (H_1)^2\\
				H_4 &= -\epsilon
			\end{alignedat}
		\right.;
		\end{equation}
		}
	 	\item{cohomogeneity $3$ if and only if 
			\begin{equation}\label{caseA1cohomog3Cond}
				\begin{alignedat}{2}
				\epsilon(3H_2H_3 +H_3&+4 H_2)\cdot((H_1)^2 - \epsilon H_2)\\
						 &= -  (\epsilon H_2 + 3(H_1)^2)(H_4+\epsilon)\cdot(H_2 - 1) 
				\end{alignedat}		 
			\end{equation}
		and it is not of cohomogeneity $2$;
		}
		\item{cohomogeneity $4$ if and only if it is not of cohomogeneity $2$ or $3$.}
	 \end{itemize}
	
\end{theorem}

\begin{remark}\label{A1discreteSymmetryRmk}
	In \eqref{caseA1oms} and \eqref{caseA1Hs}, if we replace
		\[
			\omega^0,\bar\omega^0,\gamma,\omega^1,\ldots,\omega^4;  H_1,\ldots,H_4
		\]
	  by the following expressions
		\[
			\frac{\epsilon H_2}{H_1}\bar\omega^0, \frac{H_2}{H_1}\omega^0, -\gamma,
				-\frac{H_2}{H_1}\omega^3, -\omega^4, -\frac{\epsilon H_2}{H_1}\omega^1,
				-\omega^2;
				H_1,\frac{\epsilon (H_1)^2}{H_2}, \epsilon H_4, \epsilon H_3,
		\]
	the resulting equations still hold. In particular, after the replacement, 
	the two cases in \eqref{caseA1cohomog2Cond} are swapped, 
	and \eqref{caseA1cohomog3Cond} remains the same.
	This discrete symmetry arises from switching the two underlying Monge-Amp\`ere systems as well as the 
	two characteristic systems of $(N,\B)$. Because of this, we regard the two cases contained in
	\eqref{caseA1cohomog2Cond} as equivalent.		
\end{remark}

\subsection{Coordinate forms}  Observe that, in \eqref{caseA1oms}, both $\omega^1$ and 
$\omega^3$ are integrable. Using Lemma \ref{GenAss123CharSysLemma} and 
Proposition \ref{genFGordonProp}, we immediately deduce 
			the following result.
			
	\begin{Prop}
		For any Type $\mathscr{A}_1$ B\"acklund transformation whose relative invariants satisfy
		\eqref{caseA1ThmAssu}, the underlying hyperbolic Monge-Amp\`ere systems must
		be contact equivalent to PDEs of the form
			\[
				z_{xy} = F(x,y,z,z_x,z_y);
			\]
		furthermore, one can choose coordinates such that the $x,y$ variables are 
		preserved under the B\"acklund transformation.
	\end{Prop} 

	Now suppose that a Type $\mathscr{A}_1$ B\"acklund transformation satisfies not only \eqref{caseA1ThmAssu}
	but also the `cohomogeneity $2$' condition
			\begin{equation}
				\left\{
				\begin{alignedat}{2}
					H_2 &= 1,\\
					H_3 &= -1.
				\end{alignedat}
				\right.
			\end{equation}
	By \eqref{caseA1Hs}, the zero locus of $H_4$ in $N$ has empty interior. 
	Wherever $H_4\ne 0$, we introduce
	 $1$-forms $\sigma^0,\sigma^1,\ldots, \sigma^4$ as follows.
		\[
		(\sigma^0,\sigma^1,\ldots,\sigma^4):=
				\(\frac{\omega^0}{H_1}, \frac{\omega^1}{H_4}, 
				H_4\omega^2, \frac{\omega^3}{H_1H_4}, H_4(\omega^4+\gamma)\).	
		\]
	By \eqref{caseA1oms} and \eqref{caseA1Hs}, $(\sigma^i)$ satisfy the structure equations
		\begin{equation}\label{A1cohomog2MA1}
			\left\{
			\begin{alignedat}{1}
				\ed\sigma^0& = -\epsilon \sigma^0\W(\sigma^2-\sigma^4) 
							+ \sigma^1\W\sigma^2+\sigma^3\W\sigma^4,\\
				\ed\sigma^1&= \phantom{+}\epsilon\sigma^1\W(\sigma^2+\sigma^3+\sigma^4),\\
				\ed\sigma^2&=\phantom{+}\sigma^0\W(\sigma^2-\sigma^4) 
							- \epsilon(\sigma^1-\sigma^3)\W\sigma^2,\\
				\ed\sigma^3&=-\epsilon \sigma^3\W(\sigma^1+\sigma^2+\sigma^4),\\
				\ed\sigma^4&=-\sigma^0\W(\sigma^2-\sigma^4) - \epsilon(\sigma^1-\sigma^3)\W\sigma^4.
			\end{alignedat}
			\right.
		\end{equation}
	It follows that one of the hyperbolic Monge-Amp\`ere systems being B\"acklund-related is homogeneous and has the 
	differential ideal
		\begin{equation}\label{A1cohomog2MA1ideal}
			\I = \<\sigma^0,\sigma^1\W\sigma^2,\sigma^3\W\sigma^4\>.
		\end{equation}
	
	To determine the other underlying Monge-Amp\`ere system, we introduce functions $R,S$  and 
	$1$-forms $\tau^0,\tau^1,\ldots,\tau^4$ as follows, wherever $H_4$ and $3H_4+4\epsilon$ are both nonzero.
		\begin{equation*}
			\left\{
				\begin{alignedat}{1}
				R&:= \epsilon (H_1)^2,\\
				S&:= \frac{H_4}{3H_4+4\epsilon},
				\end{alignedat}
			\right.
		\end{equation*}	
		\[			
				(\tau^0,\tau^1,\ldots, \tau^4):=\(\bar\omega^0, \frac{\epsilon \omega^1}{H_4}, 
						\epsilon H_4(\omega^2+\gamma) - \frac{S\bar\omega^0}{R},
						\frac{\epsilon \omega^3}{H_1H_4}, \epsilon H_4\omega^4 +\frac{S\bar\omega^0}{R}	
					\).
		\]
	By \eqref{caseA1oms} and \eqref{caseA1Hs}, $(\tau^i)$ and $(R,S)$ satisfy
		\begin{equation}\label{A1cohomog2MA2}
			\left\{
				\begin{alignedat}{1}
					\ed \tau^0&= \tau^0\W\(-\frac{S\tau^1}{R} 
								+ \frac{\tau^2}{S} +S\tau^3 +\tau^4\) +\tau^1\W\tau^2
								+R\tau^3\W\tau^4,\\
					\ed\tau^1&=\tau^1\W(\tau^2+\tau^3+\tau^4),\\
					\ed\tau^2&= -\frac{R+S}{R}\tau^1\W\tau^2 - \tau^2\W\tau^3 - S\tau^3\W\tau^4,\\
					\ed\tau^3&= (\tau^1+\tau^2+\tau^4)\W\tau^3,\\
					\ed\tau^4&= \frac{S}{R}\tau^1\W\tau^2 - \tau^1\W\tau^4 +(S+1)\tau^3\W\tau^4.
				\end{alignedat}
			\right.
		\end{equation}
	\begin{equation}\label{A1cohomog2MA2RS}
			\left\{
				\begin{alignedat}{1}
					\ed R&= -\frac{R(S+1)}{S}\tau^2 - \frac{R(R+S)}{S}\tau^4,\\
					\ed S&= -\frac{(R+S)S}{R}\tau^1 -(S+1)\tau^2+S(S+1)\tau^3.	
				\end{alignedat}
			\right.
		\end{equation}
	The corresponding Monge-Amp\`ere system has the differential ideal
		\begin{equation}\label{A1cohomog2MA2ideal}
			\bar \I = \<\tau^0,\tau^1\W\tau^2,\tau^3\W\tau^4\>.
		\end{equation}
	
	\begin{remark}\label{RSgradientRemark}
		Since $H_1>0$ and $\epsilon(H_1)^2\ne 1$, we have that $R\ne 0,1$.
		By \eqref{A1cohomog2MA2RS}, both $R$ and $S$ have nonzero gradients everywhere.
	\end{remark}
		
	\begin{Prop}\label{TypeA1Cohomog2MA1Coord}
		Up to contact equivalence, the hyperbolic Monge-Amp\`ere
		system characterized by \eqref{A1cohomog2MA1} and \eqref{A1cohomog2MA1ideal} 
		corresponds to the PDE:
			\begin{equation}\label{TypeA1Cohomog2MA1PDE}
				(x+y)z_{xy} + 2\sqrt{z_xz_y} =0 .	
			\end{equation}
	\end{Prop}
	
	\emph{Proof}.	 In \eqref{A1cohomog2MA1}, the two cases that correspond to $\epsilon = 1$ or $-1$
	are turned into each other if we swap $(\sigma^1,\sigma^3)$ with $(\sigma^2,\sigma^4)$, which 
	does not affect the differential ideal $\I$. Thus, we proceed assuming $\epsilon = 1$.
	
	Since $\sigma^1,\sigma^3$ are integrable, locally there exist functions $f,g,x,y$ such that
		\begin{equation}\label{A1cohomog2eta1eta3}
			\left\{
			\begin{alignedat}{1}
				\sigma^1 &= e^f\ed x, \\
				 \sigma^3 &= e^g\ed y.
			\end{alignedat}
			\right. 
		\end{equation}
	The equations of $\ed\sigma^1,\ed\sigma^3$ then imply that there exists a function $F(x,y)$ such that 
		\begin{equation}
			\left\{
			\begin{alignedat}{1}	
				f+g &= F(x,y), \\
				\sigma^2+\sigma^4 &= -\ed f + F_x\ed x - \sigma^1 - \sigma^3.
			\end{alignedat}
			\right.	
		\end{equation}
	Substituting these in the equation of $\ed(\sigma^2+\sigma^4)$, we find that $F$ must satisfy
		\begin{equation}
			F_{xy} = -2e^F.
		\end{equation}

	By using the flexibility of choosing $x,y$ in \eqref{A1cohomog2eta1eta3},
	we can eliminate the ambiguity of $F$, normalizing it to
		\begin{equation}
			F(x,y) = -2\ln|x-y|\qquad (x\ne y).
		\end{equation}
	Technically, this amounts to showing that for any $F$ satisfying \eqref{A1cohomog2eta1eta3},
	one can always find $u(x), v(y)$ that satisfy
		\[
			e^{F(x,y)} = \frac{u'(x)v'(y)}{(u(x)-v(y))^2},
		\]
	which is a consequence of the Frobenius theorem.

	To proceed, note that $\sigma^2-\sigma^4$ is integrable, so there exist functions $h, s$ such that 
		\begin{equation}\label{A1cohomog2eta24}
			\sigma^2 - \sigma^4 = e^h \ed s.
		\end{equation}
	By expanding $\ed(\sigma^2-\sigma^4)$ using \eqref{A1cohomog2MA1}, we deduce that 
		\begin{equation}
			\sigma^0 = \frac{1}{2}(\ed h + e^f \ed x - e^g\ed y + r\ed s)
		\end{equation}
	for some function $r$. 
	
	Finally, the equation of $\ed\sigma^0$ enforces that 
		\begin{equation}
			r = -e^h+G(s),
		\end{equation}
	for some function $G(s)$. 
	Using the flexibility of choosing $s$ and $h$ in \eqref{A1cohomog2eta24}, we can arrange that 
		\[
			G(s)=0.
		\]		
	
	Introducing the new variables
		\[
			\left\{ 
			\begin{alignedat}{1}
				z&:= s+ e^{-h},\\
				p& := e^{f-h},\\
				q&:= -e^{F - f -h},
			\end{alignedat}
			\right.	 
		\]		
	we obtain:	
		\begin{equation}
			\left\{ 
			\begin{alignedat}{1}
				-2e^{-h}\sigma^0 &= \ed z - p\ed x - q\ed y, \\
				2e^{-(f+h)}\sigma^1\W\sigma^2& = \sqrt{-\frac{q(x-y)^2}{p}}\ed p \W\ed x 
							+ \ed x\W\ed z+q\ed x \W\ed y.
			\end{alignedat}
			\right.
		\end{equation}			
	It follows from the vanishing of $\sigma^0$ and $\sigma^1\W\sigma^2$ that $z = z(x,y)$ satisfies
		\[
			(x-y)z_{xy} = \pm 2\sqrt{-z_xz_y}.
		\]	
	By flipping the sign of either $x$ or $y$, we obtain \eqref{TypeA1Cohomog2MA1PDE}. \qed

	 \begin{Prop}\label{TypeA1Cohomog2MA2Coord}
		Up to contact equivalence, the hyperbolic Monge-Amp\`ere
		system characterized by \eqref{A1cohomog2MA2}-\eqref{A1cohomog2MA2ideal} 
		corresponds to the PDE:
			\begin{equation}\label{TypeA1Cohomog2MA2PDE}
				(x+y)z_{xy} - \mathcal{W}(z_x,z_y) = 0,
			\end{equation}
		where \[\mathcal{W}(p,q) = ({W_0}(e^{p})+1)({W_{-1}}(-e^q)+1)\] with 
		$W_0, W_{-1}$ being the two 
		real branches of the Lambert $W$ function.
	\end{Prop}
	 
	 \emph{Proof.} The argument is similar to the proof of Proposition \ref{TypeA1Cohomog2MA1Coord},
	 			but the calculations are more complicated.
				
				First, we have chosen $(\tau^i)$ in a way so that
				the coefficient of $\tau^0$ in $\ed\tau^0$ is a closed $1$-form; that is, locally
				there exists a function $r$ such that
				\[
					\ed r = -\frac{S\tau^1}{R}+\frac{\tau^2}{S} +S\tau^3+\tau^4.
				\]
				It follows that 
				\[
					\ed(e^r\tau^0) = e^r(\tau^1\W\tau^2+ R\tau^3\W\tau^4).
				\]
				
				Second, similar to 
				how the proof of Proposition \ref{TypeA1Cohomog2MA1Coord} began, 
				we observe that $\tau^1,\tau^3, \tau^2+\tau^4$ are integrable. From this we deduce that,
				for certain independent functions $f,x,y$, the following hold:
				\begin{equation}
					\left\{
					\begin{alignedat}{1}
						\tau^1 &= e^f \ed x,\\
						\tau^3 & = e^{F-f}\ed y,\\
						\tau^2+\tau^4& = -\ed f - e^{F-f}\ed y +(F_x - e^f)\ed x,
					\end{alignedat}
					\right.
				\end{equation}
	 			where $F = F(x,y)$ must satisfy $F_{xy} = -2e^F$. Using the flexibility
				of choosing $x,y$, we can arrange that
				\begin{equation}\label{TypeA1Cohomog2MA2F}
					F(x,y)  = -2\ln|x-y|\quad  (x\ne y).
				\end{equation}	 
			
				By Remark \ref{RSgradientRemark} and by shrinking the domain, if needed, 
				we can assume that $S \ne\pm 1$.
				With this assumption, $\tau^1,\ldots,\tau^4$ are expressible as linear combinations
				of $\ed x, \ed y, \ed r, \ed f$ with coefficients being functions of $x,y,f,R,S$.
				
				Next, using the expression of $\ed S$, we find that
				\begin{equation}\label{2lnSetc}
					\ed(2\ln|S+1| - \ln|S| - f - r) = \(\frac{2e^f}{S+1}+\frac{2Se^f}{R(S+1)} +\frac{2}{x-y}\)\ed x;
				\end{equation}
	 			thus, \[2\ln|S+1| - \ln|S| - f - r = T(x)\]
				for some function $T(x)$.
				
				Let $T'$ denote coefficient of $\ed x$ on the right-hand-side of \eqref{2lnSetc}. 
				We observe by direct calculation that 
					\[
						\ed T' = -\frac{1}{2} (T')^2 \ed x.
					\]
	 			Therefore, either 
					\[
						T(x) = 2\ln|x+a|
					\]	
				or
					\[
						T(x) = b,
					\]	
	 			where $a,b$ are constants. In the latter case, we translate $r$ by a constant
				to set $T = 0$. In the former case, we first add an appropriate constant to both  $x$ and $y$ 
				 to set $a=0$, and then replace $(x,y)$ 
				by $(-x^{-1}, -y^{-1})$ to arrange that $T = 0$. In both cases, 
				\eqref{TypeA1Cohomog2MA2F} is preserved.
								
				Now $T = 0$ implies that 
					\begin{equation}\label{eqforS}
						(S+1)^2= |S|e^{f+r},
					\end{equation}
				which involves two cases: $S>0$ and $S<0$. 
								
				When $S>0$, it is necessary that $e^{f+r}>4$, or \eqref{eqforS} cannot be solved for $S$ (noting 
				that $S\ne 1$ as well). Thus, introduce a variable $\alpha$ and write
					\begin{equation}
						e^{f+r} = 2\cosh \alpha +2 \quad (\alpha >0),
					\end{equation}
				 which implies that
					\[
						S = \cosh\alpha \pm \sinh \alpha.	
					\]
	 			By allowing $\alpha$ to be negative, we have
					\begin{equation}
						S = \cosh\alpha + \sinh\alpha\quad (\alpha \ne 0).
					\end{equation}
				The equation $T'=0$ determines $R$, that is,
					\begin{equation}
						R = -\frac{(x-y)S e^f}{(x-y)e^f+S+1}.
					\end{equation}
				The equations in \eqref{A1cohomog2MA2RS}
				are now identities.	
				
				In the last steps, we express
					\[
						\ed(e^r\tau^0) =e^r(\tau^1\W\tau^2+R\tau^3\W\tau^4)
					\]
				in the form $\ed x\W\ed p + \ed y \W\ed q$, where $p = - (\alpha+e^\alpha)$, and
					\[
						q = \ln\(\frac{e^{r+\alpha}}{(e^\alpha +1)(x-y)} + 1\) - \frac{e^{r+\alpha}}{(e^\alpha +1)(x-y)} - 1.
					\] 
				Meanwhile,
					\begin{equation}\label{TypeA1Cohomog2MA2tau2Wtau2}
						e^r(\tau^1\W\tau^2) 
							= -\frac{(e^\alpha - 1)\sinh \alpha}{\cosh \alpha -1}\ed x\W\ed \alpha
							+\frac{e^{r+\alpha}}{(x-y)^2}\ed x\W\ed y,
					\end{equation}
				and
					\begin{equation}\label{TypeA1Cohomog2MA2dxWdp}
						\ed x\W\ed p = \frac{2\ed x\W\ed \alpha}{\tanh(\alpha/2)-1}.			
					\end{equation}
				Since $\tau^1\W\tau^2$ vanishes on integral manifolds, it follows 
				from \eqref{TypeA1Cohomog2MA2tau2Wtau2} and \eqref{TypeA1Cohomog2MA2dxWdp}
				that the underlying PDE is
				contact equivalent to
					\[
						z_{xy} = -\frac{e^{r+\alpha}}{(x-y)^2}.
					\]	
				By expressing $e^{r+\alpha}$ in terms of $x,y,p,q$ and by flipping the sign of $x$,
				the PDE can be written as (assuming $x+y<0$ and thus $q<-1$)
					\begin{equation}\label{TypeA1Cohomog2MA2CoordEq}
						(x+y)z_{xy} - ({W_0}(e^{p})+1)({W_{-1}}(-e^q)+1) = 0,
					\end{equation}
				where ${W_0}$ and $W_{-1}$ stand for the two real branches of the 
				\emph{Lambert $W$ function}. 
				
				When $S<0$, the argument almost identical to the above, from which we obtain 
				\eqref{TypeA1Cohomog2MA2CoordEq} with $x+y>0$.
				
				This completes the proof.\qed\\
					 	
		To sum up, we have proved the following result.
		
		\begin{theorem}\label{TypeA1Cohomog2Thm}
			If a Type $\mathscr{A}_1$ rank $2$ B\"acklund transformation whose relative invariants satisfy
			\eqref{caseA1ThmAssu} attains the maximum possible symmetry (i.e., of cohomogneity $2$), then, up to 
			contact equivalence, it relates solutions of
				\begin{equation}\label{GoursatI_inthm}
					(x+y)z_{xy} + 2\sqrt{z_xz_y} =0
				\end{equation}
			with those of the following equation
				\begin{equation}\label{LambertW_inthm}
					(x+y)z_{xy} - \mathcal{W}(z_x,z_y) = 0,
				\end{equation}
		where \[\mathcal{W}(p,q) = ({W_0}(e^{p})+1)({W_{-1}}(-e^q)+1)\] with 
		$W_0, W_{-1}$ being the two 
		real branches of the Lambert $W$ function.
		\end{theorem}					
	 \begin{remark}
	 	It is worth noting that \eqref{GoursatI_inthm} belongs to Goursat's list \cite{Goursat1899} 
		of Darboux integrable PDEs of the form $z_{xy} = F(x,y,z,z_x,z_y)$ (see also \cite{CI05}). 
		It was observed in \cite{Zvyagin1991} that there exists a rank $1$ B\"acklund transformation 
		relating \eqref{GoursatI_inthm} and the homogeneous wave equation $u_{xy} = 0$.	 
	\end{remark}

\section{Type $\mathscr{B}$}

			Let $(N,\B;\pi, \bar \pi)$ be a Type $\mathscr{B}$ rank $2$ B\"acklund transformation relating two
			hyperbolic Monge-Amp\`ere systems  $(M,\I)$ and $(\bar M,\bar \I)$. 
			Let $\chi_{10}$ and $\chi_{01}$ be an ordered pair of characteristic systems of $(N,\B)$. 
			By the assumption, for each $p\in N$, there exists an open neighborhood
			$U$ of $p$, a coframing
				\[
					\bs\omega = (\omega^0,\bar\omega^0,\gamma,\omega^1,\ldots,\omega^4)
						\in\mathcal{F}^*(U)
				\]
			and two functions $A_1,A_2$ defined on $U$ 	
			such that
			\begin{equation}\label{caseBstreqBegin}
				\left\{
				\begin{alignedat}{3}
				\ed\omega^0&\equiv \omega^1\W\omega^2+\omega^3\W\omega^4 
							&&\mod\omega^0,\\
				\ed\bar\omega^0&\equiv \omega^1\W\omega^2+\omega^3\W\omega^4 
							&&\mod\omega^0,\bar\omega^0,\gamma,\\
				\ed\gamma&\equiv A_1\omega^1\W\omega^2+A_2\omega^3\W\omega^4 
							&&\mod\omega^0,\bar\omega^0,\gamma.
				\end{alignedat}
				\right.
			\end{equation}
			The generic case is when $A_1\ne A_2$, which corresponds to $\B$ being \emph{Pfaffian}, 
			which means that $\B$ is differentially generated by its degree $1$ part. 
			This is what we will assume from now on.
			
			To `refine' the coframing, we take the following steps successively.
		
		\begin{enumerate}[\it Step 1.]
		
			\item{swap the pairs $(\omega^1,\omega^2)$ and $(\omega^3,\omega^4)$, 
				if needed, to arrange that
				\begin{equation}\label{caseBcharSysMatch}
					\chi_{10} = \lbb\omega^0,\bar\omega^0,\gamma,\omega^1,\omega^2\rbb,\qquad
					\chi_{01} = \lbb\omega^0,\bar\omega^0,\gamma,\omega^3,\omega^4\rbb.
				\end{equation}}
			\item{add a multiple of $\omega^0$ into $\gamma$ to arrange that \[A_1 = -A_2\ne 0.\]}
			\item{scale $\gamma$ to arrange that $A_1 = -A_2 = 1$; as a result,
				\begin{equation}\label{GA1gam}
					\ed \gamma \equiv \omega^1\W\omega^2 - \omega^3\W\omega^4 
						\mod \omega^0,\bar\omega^0,\gamma.
				\end{equation}}
			\item{add suitable multiples of $\omega^0$ into $\omega^1,\ldots,\omega^4$ to arrange that
				\begin{equation*}
					\ed\bar\omega^0\equiv \omega^1\W\omega^2
								+\omega^3\W\omega^4\mod\bar\omega^0,\gamma.
				\end{equation*}}
			\item{
				add suitable multiples of $\gamma$ into $\omega^3,\omega^4$ such that the 
				following congruences hold: 
				\begin{align}
					\ed\omega^0&\equiv \omega^1\W\omega^2+\omega^3\W\omega^4
						+(B_3 \omega^3+B_4\omega^4)\W\gamma &&\mod\omega^0,  \label{GA1om0}
								\\ 
					\ed\bar\omega^0&\equiv \omega^1\W\omega^2+\omega^3\W\omega^4
						+(B_1\omega^1+B_2\omega^2)\W\gamma &&\mod\bar\omega^0, \label{GA1bom0}
				\end{align}
				where $B_1,\ldots,B_4$ are functions defined on $U$. 
				(Note that $\ed\bar\omega^0$ cannot have an 
				$\omega^0\W\gamma$ term, since $\bar\omega^0\W(\ed\bar\omega^0)^3=0$.)}
		\end{enumerate}

	\begin{Def}\label{caseB0adaptedDef}					
	 Let $(N,\B;\pi,\bar\pi)$ 
		be a Type $\mathscr{B}$ rank $2$ B\"acklund transformation with an ordered pair of characteristic systems
		$\chi_{10},\chi_{01}$ for $(N,\B)$, where $\B$ is Pfaffian.
		A coframing 
		\[
			\bs\omega:=(\omega^0,\bar\omega^0,\gamma,\omega^1,\ldots,\omega^4)
		\]	 
		defined on an open subset
		 $U\subset N$ is said to be \emph{$0$-adapted}
		 if it satisfies 
 		\begin{equation}\label{caseBrank2zeroAdBundles}
 				\lbb\omega^0\rbb = \pi^*I^1, \quad 
				\lbb\bar\omega^0\rbb = \bar \pi^*\bar I^1, \quad 
				\lbb\omega^0,\bar\omega^0,\gamma\rbb = B^1,	
		\end{equation}
		the condition \eqref{caseBcharSysMatch} and 
		the congruences \eqref{GA1gam}, \eqref{GA1om0}, \eqref{GA1bom0} for some functions
		$B_1,\ldots, B_4$ defined on $U$.
\end{Def}

	Given a Type $\mathscr{B}$ rank $2$ B\"acklund transformation with a 
	$0$-adapted coframing $\bs\omega$ defined on an open subset $U\subset N$, there are two types of transformations
	\[	
		\bs\omega\mapsto\bs\omega\cdot g  = g^{-1}\bs\omega
	\]	 
	that can yield a coframing that is still $0$-adapted. One is when
	\[
		g: U\rightarrow\GL(7,\R)
	\]	 
	takes the form
		\begin{equation}\label{caseBfirstTrans}
			\begin{alignedat}{2}	
					g =\left( \begin{array}{ccc} rI_3 & 0&0\\
									        0&\bs a& 0\\
									         0& 0&\bs b
							\end{array}		       
									       \right), 
			\end{alignedat}	\qquad
			\begin{alignedat}{2}
					& \bs a = (a_{ij}), \bs b = (b_{ij}): U\rightarrow \GL(2,\R),\\
						&\qquad \det(\bs a) = \det(\bs b) = r\ne 0;
			\end{alignedat}	
		\end{equation}	
   	the other is when 
		\begin{equation}\label{caseBsecondTrans}
			{\small
			\begin{alignedat}{2}	
					g =\left( \begin{array}{ccc|cccc}
									1&0&0&0&0&0&0\\
									0&1&0&0&0&0&0\\
									c&-c&1&0&0&0&0	    \\
									\hline
									B_2c&0&0&1&0&0&0\\
									-B_1c&0&0&0&1&0&0\\
									0&-B_4c&0&0&0&1&0\\
									0&B_3c&0&0&0&0&1
							\end{array}		       
									       \right), 
			\end{alignedat}	\qquad c: U\rightarrow\R.}
		\end{equation}
 	
	All $0$-adapted coframings that are related to $\bs\omega$ by a transformation of the first type
	are precisely the local sections of a principal bundle $\mathcal{H}$ over $U$; along the fibres of $\mathcal{H}$,
	the functions $B_i$ $(i = 1,\ldots, 4)$ transform by
	\begin{equation*}
				\begin{alignedat}{2}
						\left(
							\begin{array}{c}
								B_1\\
								B_2
							\end{array}
						\right)(u\cdot g) &= \bs a^T
								\left(
									\begin{array}{c}
										B_1\\
										B_2
									\end{array}
								\right)(u),	\\	
						\left(
							\begin{array}{c}
								B_3\\
								B_4
							\end{array}
						\right)(u\cdot g) &= 
										\bs b^T
								\left(
									\begin{array}{c}
										B_3\\
										B_4
									\end{array}
								\right)(u).		
					\end{alignedat}			
				\end{equation*}
	On the other hand, \eqref{caseBsecondTrans} depends on the choice of $\bs\omega$ and 
	the functions $B_i$ $(i = 1,\ldots, 4)$ associated to it; and the transformation preserves 
	$B_i$. Thus, we encounter an issue of not being able to express \emph{all} $0$-adapted
	coframings as local sections of a single $H$-structure for some subgroup $H\subset \GL(7,\R)$. 
	This issue can be resolved by a so-called `structure reduction', which  normalizes $(B_i)$ by 
	applying transformations of the first type above.

	Indeed, given any $0$-adapted coframing $\bs \omega$ in the sense of 
	Definition \ref{caseB0adaptedDef}, we 
	can apply a transformation $\bs\omega\mapsto\bs\omega\cdot g$, where $g$ is of the form \eqref{caseBfirstTrans},
	to reduce to one of the following $4$ cases.
		\begin{enumerate}[\qquad\qquad]
			\item[$\mathscr{B}_1$:]{$B_1 = B_3 = 1$, $B_2 = B_4 = 0$;}
			\item[$\mathscr{B}_2$:]{$B_i = 0$ $(i = 1,\ldots, 4)$;}
			\item[$\mathscr{B}_3$:]{$B_2, B_3, B_4 = 0$, $B_1 = 1$;}
			\item[$\mathscr{B}_4$:]{$B_1, B_2, B_4 = 0$, $B_3 = 1$.}	
		\end{enumerate}							
	
	The cases $\mathscr{B}_3$ and $\mathscr{B}_4$ are essentially equivalent, since they 
	turn into each other as we swap the two underlying Monge-Amp\`ere systems.

\subsection{The Case $\mathscr{B}_1$}

	This is the case when, in \eqref{GA1om0} and \eqref{GA1bom0}, both $(B_1, B_2)$ and
	$(B_3, B_4)$ are nonzero, a generic subcase of Type $\mathscr{B}$.
	
	\begin{Def}\label{caseB1oneadaptDef}
		A local coframing $\bs\omega$ is said to be \emph{$1$-adapted} to
		a Type $\mathscr{B}_1$ rank $2$ B\"acklund transformation 
		if it is $0$-adapted with 
			\begin{equation}\label{caseB1Bs}
				(B_1,B_2) = (B_3, B_4) = (1,0).
			\end{equation}
	\end{Def}

	\begin{lemma}\label{caseBGstrlemma}			
	Let $(N,\B;\pi,\bar\pi)$ be a Type $\mathscr{B}_1$ rank $2$ B\"acklund transformation 
	relating two hyperbolic Monge-Amp\`ere systems $(M,\I)$ and $(\bar M,\bar \I)$,
	where $\B$ is Pfaffian.
	Each point $p\in N$ has an open neighborhood $U$, on which the $1$-adapted coframings 
	(Definition \ref{caseB1oneadaptDef}) 
	  are precisely the sections of a $G$-structure $\mathcal{G}$ on $N$, where 
	  $G\subset\GL(7,\R)$ is the subgroup formed by matrices of the form
	  	\begin{equation}\label{caseBGelements}
			{\small
			\begin{alignedat}{2}	
					g =\left( \begin{array}{ccc|cccc} 
								r&0&0&0&0&0&0\\
								0&r&0&0&0&0&0\\
								c&-c&r&0&0&0&0\\\hline
								0&0&0&1&0&0&0\\
								-c&0&0&a&r&0&0\\
								0&0&0&0&0&1&0\\
								0&c&0&0&0&b&r
							\end{array}		       
									       \right),
			\end{alignedat}	\qquad
					a,b,c\in\R;~ r\ne 0.}
		\end{equation}	
	\end{lemma}	
	{\it Proof.} Let $\bs\omega$ be a $1$-adapted coframing.
			For any function $g:U\rightarrow G$, the coframing $\bs\omega\cdot g = g^{-1}\bs\omega$ 
			is also $1$-adapted.
			
			Conversely, suppose that $\bs\omega$ and $\bs\omega\cdot g$ are both
			$1$-adapted coframings defined on $U$ with $g: U\rightarrow \GL(7,\R)$.
			Using \eqref{caseBcharSysMatch}-\eqref{caseBrank2zeroAdBundles} and \eqref{caseB1Bs},
			it is straightforward to check that $g$ must be of the form \eqref{caseBGelements}.
	 \qed \\

			 Let $\bs\omega$ also denote the tautological $1$-form on $\mathcal{G}$. 
		 The following structure equations hold.
			\begin{equation}\label{caseB1streqn}
				{\small
				\ed\left(\begin{array}{c}
						\omega^0\\
						\bar\omega^0\\
						\gamma\\
						\omega^1\\
						\omega^2\\
						\omega^3\\
						\omega^4
					\end{array}
				\right)=
								-\left(
									\begin{array}{ccccccc}
										\phi &0&0&0&0&0&0\\	
										0&\phi&0&0&0&0&0\\	
										\varpi&-\varpi&\phi&0&0&0&0\\
										0&0&0&0& 0&0&0\\	
										-\varpi&0&0&\alpha& \phi&0&0\\	
										0&0&0&0&0&0&0\\	
										0&\varpi&0&0&0&\beta&\phi
									\end{array}
								\right)\W
				\left(\begin{array}{c}
						\omega^0\\
						\bar\omega^0\\
						\gamma\\
						\omega^1\\
						\omega^2\\
						\omega^3\\
						\omega^4
					\end{array}
				\right) 
				+
				\left(
					\begin{array}{c}
						\Omega^0\\
						\bar\Omega^0\\
						\Gamma\\
						\Omega^1\\
						\Omega^2\\
						\Omega^3\\
						\Omega^4
					\end{array}
				\right).}
			\end{equation}
			Here, by adding a linear combination of the semibasic $1$-forms to $\phi, \varpi$ and by the 
			reproducing property of $\bs\omega$, we can arrange
			that
				\begin{equation}\label{caseB1Tor}\left\{
					\begin{alignedat}{2}
					\Omega^0&= \omega^1\W\omega^2 +\omega^3\W(\omega^4+\gamma) +(P_0 \bar\omega^0+K\gamma+P_i\omega^i)\W\omega^0,\\
					\bar\Omega^0 & = \omega^1\W(\omega^2+\gamma)+\omega^3\W\omega^4+(P_0 \omega^0 - K\gamma+ Q_i\omega^i)\W\bar\omega^0,\\
					\Gamma & = 
						\omega^1\W\omega^2-\omega^3\W\omega^4 +C_i\omega^i\W(\omega^0+\bar\omega^0),
				\end{alignedat}\right.
				\end{equation}
			for certain functions $C_i, P_i, Q_i, P_0, K$ $(i = 1,\ldots, 4)$ defined on $\mathcal{G}$.

			It is not difficult to see that the conclusions of Lemma \ref{GenAss123CharSysLemma} and
			Lemma \ref{Rank2GenericInvLemma}
			apply to the current case without any change. As a result, 
			$\Omega^1,\ldots,\Omega^4$ are the same as  
			\eqref{caseA1Om1to4} except that each $D_i$ are now replaced with $C_i$. 
			Moreover, $T^k_{ij} = -T^k_{ji}$
			for all $i,j,k = 1,\ldots,4$, 
			and 
			\[
			T^1_{23}, T^1_{24},T^1_{34}, T^2_{23}, T^2_{24},T^2_{34},T^3_{12}, T^3_{14}, T^3_{24}, 
			T^4_{12}, T^4_{14}, T^4_{24}
			\] 
			are all zero.
			
			Furthermore, one can add semibasic $1$-forms into $\alpha$ to arrange that
						\[
							  T^2_{12} = T^2_{13} = T^2_{14} = T^2_1= \bar T^2_1 =  R^2_1=0;
						\]	
			similarly, by adjusting $\beta$, we can arrange that
						\[
							  T^4_{13} = T^4_{23} = T^4_{34} = \bar T^4_3 =  T^4_3 =  R^4_3 = 0.
						\]	
			By adding a multiple of $\omega^0 - \bar\omega^0$ to $\varpi$, we can arrange that
				\[
					T^4_0  = - T^2_0.
				\]			
			The torsion cannot be absorbed further.
			
			By expanding $\ed(\ed\bs\omega)$, we can uncover relations among
			the torsion functions.
			In fact, by computing
				\begin{equation*}\left\{
					\begin{alignedat}{2}
					\ed(\ed\omega^0)&\mod\omega^0, \quad \\
					\ed(\ed\bar\omega^0)&\mod\bar\omega^0,\quad\\
					\ed(\ed\gamma)&\mod\omega^0,\bar\omega^0,\gamma,
					\end{alignedat}\right.
				\end{equation*}
			we obtain the relations below.
				\begin{equation}\label{caseB1firstrel}
					{\small\left\{\;
					\begin{alignedat}{5}
						T^2_2 &=P_0 - T^1_1+C_2,\qquad
									&&Q_1 && = P_1+1,\\
						T^4_4& = P_0 - T^3_3,\qquad
									&&Q_2&& = P_2,\\
						\bar T^2_2& = P_0- \bar T^1_1,\qquad
									&&Q_3&&= P_3+1,\\
						\bar T^4_4& = P_0 - \bar T^3_3 + C_4,	\qquad
									&&Q_4&& = P_4,\\
						T^1_{12}& = -2K - P_2-C_4  +\frac{P_4}{2},\qquad
									&&S_1&& = C_2+\frac{P_2}{2},\\
						T^3_{34}& =  2K - P_4 + C_2 + \frac{P_2}{2},\qquad
									&&S_2&& = -C_1- \frac{P_1}{2},\\
						T^1_{13}& = -C_3 - C_1- \frac{P_1+P_3+1}{2},\qquad
									&&S_3&& = -C_4+\frac{P_4}{2},\\
						T^3_{13}& = -C_1 - C_3 + \frac{P_1+P_3+1}{2},\qquad
									&&S_4&& = C_3 - \frac{P_3+1}{2},\\
						T^1_{14}& = -C_4 - \frac{P_4}{2},\qquad
									&&R^2_2&& = K+C_4- R^1_1 - \frac{P_4}{2},\\
						T^3_{23}& = -C_2+\frac{P_2}{2},\qquad
									&&R^4_4&& = -K-C_2-R^3_3- \frac{P_2}{2}.
					\end{alignedat}\right.}
				\end{equation}		
			
			There are $25$  torsion functions remaining, and they are not independent.
			
			By computing
					\begin{equation*}\left\{
						\begin{alignedat}{2}
						\ed(\ed\omega^0) &\mod \bar\omega^0,\gamma,\omega^2,\omega^3,\\
						\ed(\ed\bar\omega^0)&\mod\omega^0,\gamma,\omega^2,\omega^3,
						\end{alignedat}\right.
					\end{equation*}
			we obtain two expressions for $\ed\phi$, both reduced modulo $\{\omega^0,\bar\omega^0,
			\gamma,\omega^2,\omega^3\}$:
				\begin{equation}\label{caseB1dphiexpr}
					\left\{
					\begin{alignedat}{2}
						\ed\phi&\equiv C_4\omega^1\W\omega^4 +\psi\\
						\ed\phi&\equiv -3C_4\omega^1\W\omega^4+\psi
					\end{alignedat}\right\}\mod\omega^0,\bar\omega^0,\gamma,\omega^2,\omega^3,
				\end{equation}
			where
				\begin{equation*}
					\begin{alignedat}{2}
						\psi& = (\ed P_1+\varpi - P_2\alpha)\W\omega^1 
									  + (\ed P_4 - \phi P_4)\W\omega^4\\
							& - \(\frac{P_1P_4}{2} + P_1C_4+P_3C_4
								 	 - P_4C_3+\frac{P_4}{2}\)\omega^1\W\omega^4		  .
					\end{alignedat}
				\end{equation*}
			From \eqref{caseB1dphiexpr} we deduce that
				\begin{equation}\label{caseB1C4zero}
					C_4 = 0.
				\end{equation}

			Similarly, by expanding
					\begin{equation}\label{caseB1ddom0ij}
						\left\{
						\begin{alignedat}{2}
						\ed(\ed\omega^0) &\mod \bar\omega^0,\gamma,\omega^i,\omega^j,\\
						\ed(\ed\bar\omega^0)&\mod\omega^0,\gamma,\omega^i,\omega^j,
						\end{alignedat}\right.
					\end{equation}
			for $(i,j)$ being $(1,4)$, $(2,4)$ and $(3,4)$ and by comparing in each case the two expressions
			of $\ed\phi$ (reduced modulo the same set of $1$-forms), we obtain
				\begin{equation}\label{caseB1C2C3K}
					\left\{
					\begin{alignedat}{2}
						C_2& = 0,\\
						C_3& = -C_1,\\
						K& = \frac{P_4 - P_2}{4}.
					\end{alignedat}
					\right.
				\end{equation}
			Computing \eqref{caseB1ddom0ij} for $(i,j)$ being $(1,2)$ and $(1,3)$ does not yield new relations.	
						
			Furthermore, computing $\ed(\ed\gamma)$ and reducing in two ways:
				\begin{equation}\label{caseB1ddgam}\left\{
					\begin{alignedat}{2}
						\ed(\ed\gamma)&\mod \omega^0, \gamma,\omega^1 - \omega^3,\\
						\ed(\ed\gamma)&\mod \bar\omega^0,\gamma,\omega^1 - \omega^3,
					\end{alignedat}\right.
				\end{equation}			
			we obtain two expressions for $\ed\varpi$:
				\begin{equation}\label{caseB1dvarpiexpr}
					\left\{
					\begin{alignedat}{2}
						\ed\varpi&\equiv \xi + \frac{\zeta}{2}\\
						\ed\varpi&\equiv \xi - \frac{\zeta}{2}
					\end{alignedat}\right\}\mod\omega^0,\bar\omega^0,\gamma,\omega^1 - \omega^3,
				\end{equation}
			where
				\begin{equation}\label{caseB1formofxizeta}
					\left\{
						\begin{alignedat}{2}
							\xi & = \varpi\W(P_2\omega^2+(P_1+P_3+1)\omega^3+ P_4\omega^4),\\
							\zeta& = (2P_0+ C_1(P_2 - P_4))\omega^2\W\omega^3\\
								& +(2P_0 - C_1(P_2 - P_4))\omega^3\W\omega^4.
						\end{alignedat}
					\right.
				\end{equation}
			From \eqref{caseB1formofxizeta} and \eqref{caseB1C2C3K}
			we deduce
				\begin{equation}\label{P0C1P2P4}\left\{
				\begin{alignedat}{2}
					P_0  &= 0,\\
					C_1K& = 0.
				\end{alignedat}\right.
				\end{equation}
			In particular, we have either $C_1 = 0$ or $K  = 0$ (or both), but we do not make a choice yet.
			
			We use the structure equations to determine
			 the infinitesimal transformation of the $20$ remaining torsion functions
			under the $G$-action, listed as follows.						
				
			In \eqref{caseB1torVar}, all congruences hold modulo 
			$\{\omega^0,\bar\omega^0,\gamma,\ldots,\omega^4\}$.		
			\begin{equation}\label{caseB1torVar}
				\left\{\;
				\begin{alignedat}{1}
				&	\begin{alignedat}{5}
						&\ed C_1&&\equiv 0, \\
									& \ed P_1&&\equiv \alpha P_2 - 2\varpi,\qquad\quad
									&&\ed R^1_1&&\equiv \phi R^1_1 +\alpha R^1_2,\\
									& \ed P_2 &&\equiv  \phi P_2,\quad
									&&\ed R^1_2&&\equiv 2\phi R^1_2,\\
									& \ed P_3&&\equiv  \beta P_4 + 2\varpi,\quad
									&&\ed R^3_3 &&\equiv \phi R^3_3 +\beta R^3_4,\\
									&	\ed P_4&&\equiv \phi P_4,\quad
									&&	\ed R^3_4&&\equiv 2\phi R^3_4,	
					\end{alignedat}\\
				&	\begin{alignedat}{1}
						\ed T^1_0&\equiv 2\phi T^1_0 - \varpi(T^1_2 - \bar T^1_2),\\
						\ed T^2_0&\equiv -\frac{\alpha}{2}T^1_0+\phi T^2_0 +\frac{\beta}{2}T^3_0
									+\frac{\varpi}{2}(T^1_1 - T^3_3 - \bar T^1_1+\bar T^3_3),\\
						\ed T^3_0&\equiv 2\phi T^3_0 - \varpi(T^3_4 - \bar T^3_4),		\\
						\ed T^1_1&\equiv \phi T^1_1 + \alpha T^1_2 
									-\varpi\(\frac{P_2}{2}- R^1_1\),\\
						\ed T^1_2&\equiv 2\phi T^1_2+\varpi R^1_2,\\
						\ed \bar T^1_1&\equiv \phi \bar T^1_1 +\alpha \bar T^1_2
									+\varpi\(\frac{P_4}{2} - R^1_1\),\\
						\ed\bar  T^1_2&\equiv 2\phi \bar T^1_2 - \varpi R^1_2,\\
						\ed \bar T^3_3&\equiv \phi \bar T^3_3 +\beta \bar T^3_4
									 - \varpi\(- \frac{P_4}{2} + R^3_3\),\\
						\ed\bar T^3_4&\equiv 2\phi \bar T^3_4 - \varpi R^3_4,\\
						\ed T^3_3&\equiv \phi T^3_3+\beta T^3_4+\varpi\(- \frac{P_2}{2} + R^3_3\),\\
						\ed T^3_4&\equiv 2\phi T^3_4 +\varpi R^3_4.	
					\end{alignedat}
				\end{alignedat}	
				\right.
			\end{equation}
	In other words, to each point in $N$ is associated a $20$-dimensional real representation of $G$.
	By the expressions of $\ed P_1$ and $\ed P_3$, any $G$-orbit must be at least $1$-dimensional. 
	Furthermore, we observe the following.
	
	\begin{Prop}\label{caseB12DorbitProp}
		The $G$-action
		represented by \eqref{caseB1torVar} has all its orbits being at most $2$-dimensional
		 if and only if, in \eqref{caseB1streqn}, both $\omega^1$ and $\omega^3$ are integrable.
	\end{Prop}
		
	\emph{Proof}. Consider the following sets of functions:
			\begin{equation*}
				\begin{alignedat}{1}
				\mathfrak{C}_0&=\{P_2, P_4\},\\
				\mathfrak{C}_1&=\{R^1_2, R^3_4\},\\
				\mathfrak{C}_2&=\{ T^1_2, \bar T^1_2, T^3_4,\bar T^3_4\},\\
				\mathfrak{C}_3&=\{T^1_0, T^3_0\}.
				\end{alignedat}
			\end{equation*}
		By the expressions of $\ed P_i$ $(i = 1,\ldots, 4)$, if either $P_2$ or $P_4$ is nonzero, then
		$G$ has an orbit of dimension at least $3$. If both functions in $\mathfrak{C}_0$ are zero
		with either function in $\mathfrak{C}_1$ being nonzero, say, $R^1_2\ne 0$,
		 then there exists a $G$-orbit whose projection to the $(P_1P_3R^1_1R^1_2)$-space 
		 is at least $3$-dimensional. Similar arguments work when we assume 
		 either of the following:
		 	\begin{itemize}
				\item{both functions in $\mathfrak{C}_1$ are zero, but not all functions in 
						$\mathfrak{C}_2$ are zero;}
				\item{all functions in $\mathfrak{C}_1$ and $\mathfrak{C}_2$ are zero, but a
					function in $\mathfrak{C}_3$ is nonzero.}	
			\end{itemize}
		In each of these cases, $G$ must have an orbit of dimension greater than $2$.	
		
		On the other hand, if all functions in $\mathfrak{C}_i$ $(i = 0, \ldots, 3)$ are identically zero, then
		 neither $\alpha$ nor $\beta$ occur in the right-hand-side of \eqref{caseB1torVar}.
		In this case every $G$-orbit is at most $2$-dimensional.
		
		Finally, we calculate and obtain
			\begin{equation}\left\{
				\begin{alignedat}{2}
					\ed\omega^1&\equiv T^1_0\omega^0\W\bar\omega^0 +T^1_2\omega^0\W(\gamma
							+ \omega^2) +\bar T^1_2\bar\omega^0\W\omega^2
							\\
							&+\gamma\W\(R^1_2\omega^2  - \frac{P_2}{2}\omega^3\)
							+\frac{P_2}{2}\omega^3\W(\omega^2+\omega^4)&&\mod\omega^1,\\
					\ed\omega^3&\equiv  T^3_0\bar\omega^0\W\omega^0 
							+\bar T^3_4\bar\omega^0\W(\gamma+\omega^4)
							+ T^3_4\omega^0\W\omega^4\\
							&+\gamma\W\(R^3_4\omega^4 - \frac{P_4}{2}\omega^1\)
							+\frac{P_4}{2}\omega^1\W(\omega^2+\omega^4)&&\mod\omega^3.	
				\end{alignedat}\right.
			\end{equation}
		It follows that 
			\begin{equation*}\left\{
				\begin{alignedat}{1}
					\ed\omega^1&\equiv 0 \mod \omega^1\\
					\ed\omega^3&\equiv 0\mod\omega^3
				\end{alignedat}\right.
			\end{equation*}
		holds if and only if all functions in 
		$\mathfrak{C}_i$ $(i = 0, \ldots, 3)$	are zero.
		This completes the proof.\qed
		
		\begin{remark}
			Recall how $\omega^1$ and $\omega^3$ occur in the characteristic systems of 
			$\pi^*\I$ and $\bar\pi^*\bar \I$. From this it is not difficult to see that
			if both $\omega^1$ and $\omega^3$ are integrable, then under a suitable choice of local
			coordinates on $N$
			the B\"acklund transformation is one that relates a pair PDEs of the form
			\[
				z_{xy} = F(x,y,z,z_x,z_y)\quad \text{and}\quad Z_{XY} = G(X,Y,Z,Z_X,Z_Y),
			\]
			in such a way that $x = X$ and $y = Y$ on corresponding solutions; the converse is also true.
			As one can verify, the classical $1$-parameter family of Tzitz\'eica transformations
			belong to this subclass of Type $\mathscr{B}_1$ B\"acklund transformations.
		\end{remark}

\section{Acknowledgements}
I would like to thank Professor Robert Bryant for his guidance and encouragement. 
I'm also grateful to Professor Jeanne Clelland for her support
during and after my years at CU Boulder.

\bibliographystyle{alpha}
\normalbaselines 

\newcommand{\etalchar}[1]{$^{#1}$}

\end{document}